\documentclass{article}
\usepackage{arxiv}
\pdfoutput=1
\usepackage[utf8]{inputenc} 
\usepackage[T1]{fontenc}    
\usepackage{hyperref}       
\usepackage{url}            
\usepackage{booktabs}       
\usepackage{amsfonts}       
\usepackage{nicefrac}       
\usepackage{microtype}      
\usepackage{lipsum}		
\usepackage{graphicx}
\usepackage{doi}

\usepackage{algorithm,algorithmic}
\usepackage{bbm}
\usepackage{pifont}       
\usepackage{bbding}       
\usepackage{fontawesome}  

\usepackage{subfigure}

\usepackage{amsmath,amscd,amsbsy,latexsym,bm,amsthm,amssymb}
\usepackage{wrapfig}
\usepackage{multirow}
\usepackage{cleveref}
\usepackage{verbatim}
\usepackage{dsfont}

\newtheorem{lemma}{Lemma}
\newtheorem{remark}{Remark}
\newtheorem{example}{Example}

\newcommand{\cF}{\mathcal{F}}

\newcommand{\cN}{\mathcal{N}}

\newcommand{\abs}[1]{\left| #1 \right|}

\newcommand{\cI}[1]{\mathbb{I}}

\usepackage{color}

\title{Weak Collocation Regression for Inferring Stochastic Dynamics with L\'{e}vy Noise}

\author{ 
Liya Guo\\
Yau Mathematical Sciences Center\\
 Tsinghua University\\
 Beijing, 100084, China\\
\texttt{gly22@mails.tsinghua.edu.cn} \\
	\And	
 Liwei Lu\\
Yau Mathematical Sciences Center\\
 Tsinghua University\\
 Beijing, 100084, China.\\	\texttt{llw20@mails.tsinghua.edu.cn} \\
 \AND
 Zhijun Zeng\\
Yau Mathematical Sciences Center\\
 Tsinghua University\\
 Beijing, 100084, China\\	 \texttt{zengzj22@mails.tsinghua.edu.cn}
\AND
Pipi Hu\\
Microsoft Research AI4Science\\ Beijing, 100080, China\\ \texttt{pisquare@microsoft.com}
 \AND
\thanks{Corresponding author}
Yi Zhu\\
Yau Mathematical Sciences Center,
 Tsinghua University\\
 Beijing, 100084, China\\
Yanqi Lake Beijing Institute of Mathematical Sciences and Applications\\ Beijing, 101408, China\\
\texttt{yizhu@tsinghua.edu.cn}
}

\renewcommand{\shorttitle}{\textit{arXiv} Template}

\hypersetup{
pdftitle={A template for the arxiv style},
pdfsubject={q-bio.NC, q-bio.QM},
pdfauthor={David S.~Hippocampus, Elias D.~Striatum},
pdfkeywords={First keyword, Second keyword, More},
}

\begin{document}
\maketitle


\begin{abstract}
With the rapid increase of observational, experimental and simulated data for stochastic systems, tremendous efforts have been devoted to identifying governing laws underlying the evolution of these systems. Despite the broad applications of non-Gaussian fluctuations in numerous physical phenomena, the data-driven approaches to extracting stochastic dynamics with L\'{e}vy noise are relatively few. In this work, we propose a Weak Collocation Regression (WCR) to explicitly reveal unknown stochastic dynamical systems, i.e., the Stochastic Differential Equation (SDE) with both $\alpha$-stable L\'{e}vy noise and Gaussian noise, from discrete aggregate data. This method utilizes the evolution equation of the probability distribution function, i.e., the Fokker-Planck (FP) equation. With the weak form of the FP equation, the WCR constructs a linear system of unknown parameters where all integrals are evaluated by Monte Carlo method with the observations. Then, the unknown parameters are obtained by a sparse linear regression. For a SDE with L\'{e}vy noise, the corresponding FP equation is a partial integro-differential equation (PIDE), which contains nonlocal terms, and is difficult to deal with. The weak form can avoid complicated multiple integrals. Our approach can simultaneously distinguish mixed noise types, even in multi-dimensional problems. Numerical experiments demonstrate that our method is accurate and computationally efficient.
\end{abstract}

\maketitle

\section{Introduction}
\label{sec:intro}

Learning the laws behind time-series data has been a very hot topic. These laws are often described by ordinary differential equations (ODEs) and partial differential equations (PDEs) \cite{long2018pde,yu2021onsagernet,hu2022revealing}. Depending on the physical problems, differential equations with stochastic terms are also used \cite{ma2021learning, chen2021solving, dai2020detecting,dai2023nonparametric}.
In this context, we consider the time-series data which can be properly modeled by stochastic differential equations (SDEs), in line with prior studies \cite{ma2021learning, chen2021solving,lu2024weak}. These SDEs typically encompass drift and diffusion terms, both being functions of state variables, while the diffusion terms are often considered to be arisen from Gaussian noise. Common methodologies employed to reveal unknown systems driven by Gaussian noise include parameter or non-parametric inference \cite{dai2020detecting, boninsegna2018sparse, picchini2014inference, rajabzadeh2016robust, garcia2017nonparametric, opper2019variational, dai2023nonparametric}, and neural network-based approaches \cite{kong2020sde, ma2021learning}.

Phenomena such as stock price fluctuations and abnormal diffusion \cite{emmer2004optimal, applebaum2009levy} display heavy-tailed distributions and jumps, which cannot be precisely characterized by SDEs driven solely by Gaussian noise. Consequently, there is a shift in focus to consider SDEs with L\'{e}vy noises, exhibiting stationary and independent increments that depend solely on the time interval. Recent studies have revealed the coexistence of Gaussian and L\'{e}vy transports in various complex systems \cite{bucher2015efficient, zan2021first, zan2020stochastic}. However, there is relatively little research on the inverse problem of dynamic systems containing L\'{e}vy noise.

We consider that the stochastic system follows SDEs driven by both Gaussian and L\'{e}vy noises, characterized by sufficient parameters to capture its dynamics. Challenges posed by L\'{e}vy noise, due to its jump properties which lead to an infinite computational domain, necessitate the consideration of a bounded area \cite{yanovsky2000levy, gao2016dynamical, yang2020generative}. Some previous works apply Kramers-Moyal formula to directly represent the unknown terms \cite{li2021data,li2022extracting}. Leveraging the fact that the data distribution from an SDE adheres to a specific Fokker-Planck (FP) equation \cite{schertzer2001fractional,duan2015introduction}, \cite{chen2021solving,yang2020generative} to infer the dynamic systems by neural networks. Motivated by \cite{chen2021solving}, we reveal the SDEs based on the FP equation. Therefore, estimating the unknown parameters of the SDEs is transformed into inferring the corresponding parameters of the PDEs. However, PDEs arising from dynamics with L\'{e}vy noise, known as Partial Integro-Differential Equations (PIDEs), incorporate nonlocal integral terms. They are also classified as fractional partial differential equations (FPDEs) due to the presence of fractional-order derivatives. 
Although numerical methods such as the finite difference method are employed to solve low-dimensional PIDEs \cite{gao2016dynamical}, extending these methods to high dimensions may result in the curse of dimensionality. Recent studies \cite{pang2019fpinns,guo2022monte} demonstrate the application of machine learning methods for the inversion of FPDEs in high-dimensional problems.

In many practical scenarios, trajectory data may not be readily accessible, necessitating the inference of a stochastic dynamical system from discrete observations. Models such as physics-informed neural networks (PINNs) and Wasserstein generative adversarial networks (WGANs) have been used to learn inversions of dynamical systems from discrete sparse samples \cite{ma2021learning,yang2020generative,chen2021solving}.
These methods facilitate the estimation of probability density functions at observation moments and the inference of unknown parameters or functional shapes. However, representing these unknown terms as a neural network is time-consuming and cannot yield an explicit function. Therefore, the current work aims to efficiently and explicitly learn dynamic systems with L\'{e}vy noise directly from a limited number of discrete observed samples.

Our work focuses on uncovering the laws governing stochastic dynamics with partially known or unknown terms. Building upon the weak collocation regression (WCR) \cite{lu2024weak}, we extend its scope to handle L\'{e}vy noise. Our approach involves the consideration of the weak form of the Fokker-Planck equation, employing operators, such as high-order derivatives, to act on predetermined kernel functions, as opposed to the unattainable density function. By approximating the weak formula using a discrete form of samples and employing Monte Carlo (MC) methods, revealing unknown parameters within a PDE is transformed into a linear regression problem. As mentioned before, the FP equation of the SDE with L\'{e}vy noise contains non-local terms, for which we apply a confluent hypergeometric function to reduce computational burden. Furthermore, our method offers several advantages.


\begin{itemize}
\item[(i)] Accuracy. Our approach can effectively identify two different types Gaussian and L\'{e}vy noise. The heavy-tailed L\'{e}vy process exerts a more significant influence than Gaussian noise, thereby introducing complexity to isolate the Gaussian component. In addition, unknown terms in stochastic systems may not be able to do expansions with basis selected in the dictionaries, and we give an example where the drift term of complex SDEs can be well approximated with pre-selected basis.

\item[(ii)] Efficiency. First, it requires fewer data points as the WCR method avoids an exponential increase in sample size as the dimensionality grows, while still maintaining a reasonable level of accuracy. Additionally, it ensures a shorter running time by allowing for the transformation of the learning of explicit representations of unknown parameters into solving a simple regression problem, without using traditional numerical formats or training complex neural networks. This quality makes it useful as a pre-training model for dynamical system learning.

\item[(iii)] High dimension. Our approach performs well in five-dimensional problems with laptop, while most previous works only conduct experiments with L\'{e}vy noise in at most two dimension \cite{chen2021solving, gao2016fokker}.

\end{itemize}

The organization of our study is as follows. Section \ref{sec:sde_fp} presents the FP equation associated with SDEs that incorporate both Gaussian noise and L\'{e}vy noise. Section \ref{sec:method} elucidates how our methodology achieves the estimation of unknown parameters. Following that, Section \ref{sec: exper} demonstrates the good performance of our approach with L\'{e}vy noise in both low- and high-dimensional settings. These experiments encompass both independent and coupled systems. Finally, Section \ref{sec:conclu} summarizes our primary conclusions and discussions.

\section{SDE with L\'{e}vy Noise and Its Fokker-Planck Equation}
\label{sec:sde_fp}

By the L\'{e}vy-Ito decomposition theorem \cite{applebaum2009levy}, a large class of random fluctuations $\bm{X}_t \in \mathbb{R}^d$ are indeed modeled as linear combinations of mutually independent Gaussian noise $\bm{B}_t \in \mathbb{R}^d$ and $\alpha$-stable L\'{e}vy process $\bm{L}_t$, whose form is mathematically expressed as \eqref{equa:sde_mf_bmlevy}.
\begin{equation}
  \begin{aligned}
\label{equa:sde_mf_bmlevy}
    d\bm{X}_t &= \bm{m}(\bm{X}_t) dt + \bm{\sigma}(\bm{X}_t) d \bm{B}_t + \bm{\xi}(\bm{X}_t) d \bm{L}_t\,,
\end{aligned}  
\end{equation}
where the drift term  $\bm{m}(\bm{X}_t)\in \mathbb{R}^d$, the diffusion term $\bm{\sigma}(\bm{X}_t) \in \mathbb{R}^{d\times d}$, and the noise intensity $\bm{\xi}(\bm{X}_t) \in \mathbb{R}^{d\times d}$ for L\'{e}vy noise are dependent on spatial position and are independent of $t$, which is consistent with previous work \cite{yang2020generative,li2022extracting}. The term $\bm{L}_t \in \mathbb{R}^d$ represents an $\alpha$-stable L\'{e}vy process in a non-negative Borel L\'{e}vy measure $\nu \in \mathcal{B}(\mathbb{R}^d)$, where $\alpha \in(0, 1) \cup (1, 2)$. Each L\'{e}vy process $\bm{L}_t$ is determined by a generating triplet $(D, \nu, b)$, with $\nu$ representing the L\'{e}vy measure governing $\bm{L}_t$, $D$ symbolizing the covariance matrix or diffusion matrix, and $b$ indicating the drift vector of the L\'{e}vy characteristic function.
We consider symmetric $\alpha$-stable $\bm{L}_t$ with $(0, \nu, 0)$. A L\'{e}vy process $\bm{L}_t$ in $\mathbb{R}^d$ is \textit{symmetric} if $\bm{L}_t \overset{d}{=} -\bm{L}_t$, where $\overset{d}{=}$ indicates that the variables at the left and right ends of the equal sign have the same distribution. For further details about the $\alpha$-stable L\'{e}vy process, refer to \cite{duan2015introduction}.

For an SDE with independent L\'{e}vy noise across each dimension, the L\'{e}vy measure $\nu$ above can be depicted as $\nu(d\bm{x}) = \nu(dx_1, dx_2, \cdots, dx_d) = \sum \limits_{i=1}^{d} \left(\nu(d x_i) \prod \limits_{k=1, k\neq i}^d \delta_0 (d x_k)\right)$ \cite{sun2017governing}, where $\delta_0 (x)$ is the Dirac function centered at $0$. The noise intensity of L\'{e}vy noise is presumed to be $\bm{\xi}(\bm{X}_t) = 
diag \{\xi_1, \xi_2, \cdots, \xi_d \}$, where $\xi_i$ are constants and may not be identical. 

To identify the unknown parameters $\bm{m}(\bm{X}_t), \bm{\sigma}(\bm{X}_t), \bm{\xi}(\bm{X}_t)$, we leverage the conclusion that the solution to the corresponding Fokker-Planck equation \cite{yanovsky2000levy, duan2015introduction,gao2016fokker} provides a transition probability density $p(\bm{x}, t)$ for $\bm{x}$ in a SDE \eqref{equa:sde_mf_bmlevy}.
Let the matrix function $\bm{G}(\bm{x}) = \frac{1}{2} \bm{\sigma}(\bm{x})\bm{\sigma}(\bm{x})^{T}$. 
The FP equation of SDE \eqref{equa:sde_mf_bmlevy} is given in equation \eqref{equa: fp_ind} \cite{chen2021solving,sun2017governing}.
\begin{equation}
    \begin{aligned}
    \label{equa: fp_ind}
    \frac{\partial}{\partial t} p(\bm{x}, t) 
    = -div  \left(\bm{m}(\bm{x}) p(\bm{x}, t)\right) + \sum_{i,j}^{d} \partial_{ij} \left( G_{i,j}(\bm{x}) p(\bm{x}, t) \right) - \sum_{j=1}^d |\xi_j|^{\alpha} \frac{\partial^\alpha}{\partial x_j^\alpha} p(\bm{x}, t) \,,
\end{aligned}
\end{equation}
with initial distribution $p(x, 0) = p_0(x)$, where $G_{i,j}$ is the $i$-th row and the $j$-th column of the matrix $\bm{G}$. The operator $\partial_{ij}$ is a second-order derivative operator that takes the first derivative of the $i$-th variable and the $j$-th variable respectively, and $\frac{\partial^\alpha}{\partial x_j^\alpha}$ is the $\alpha$-order fractional derivative of the $j$-th variable. The fractional Laplacian $(-\Delta)^{\alpha/2}$ of a scalar function $p$ on $\mathbb{R}^d$ is defined as equation \ref{equa: fourier_frac_lap1} \cite{burkardt2021unified,kwasnicki2017ten}. The fractional derivative $\frac{\partial^\alpha}{\partial x^\alpha} = -(-\Delta)^{\frac{\alpha}{2}} p$ in one dimension.
\begin{equation}
    \begin{gathered}
    \label{equa: fourier_frac_lap1}
-(-\Delta)^{\frac{\alpha}{2}} p(\bm{x}) = C_{d, \alpha} P.V. \int_{\mathbb{R}^d} \frac{p(\bm{x} - \bm{y})}{|\bm{x} - \bm{y}|_2^{d+\alpha}} d\bm{y} \,,
\end{gathered}
\end{equation}
where $C_{d, \alpha} = \frac{2^{\alpha-1}\alpha \Gamma(\frac{\alpha+d}{2})}{\pi^{d/2} \Gamma(1-\alpha/2)}$, P.V. means the Cauchy principal value of the integral, and $\abs{\bm{x} - \bm{y}}$ is the Euclidean distance. It has an equivalent definition that
\begin{equation}
    \begin{gathered}
    \label{equa: fourier_frac_lap}
-(-\Delta)^{\frac{\alpha}{2}} p(\bm{x}) = \cF^{-1}\left[\|\bm{k}\|^\alpha \cF[p](\bm{k}) \right] \,,
\end{gathered}
\end{equation}
where $\cF[p]$ is the Fourier transform of $p(\bm{x})$ and $\bm{k} = (k_1, \cdots, k_d)^T$ is a vector.

\section{Methodology}
\label{sec:method}

\subsection{Overview}
\label{sec:overview}
We consider the aggregate data, which refers to a data format, in which samples at certain discrete time instants are available, while the full trajectory of each individual is not available. With $L$ snapshots, $N_{t_\ell}$ samples at each snapshot $t_\ell, \ell = 1, 2, \cdots, L$ are recorded as $\mathbb{X}_{t_\ell} = \{\bm{x}_{t_\ell}^k \}_{k=1}^{N_{t_\ell}}$, and the set of the aggregate data is denoted as $\mathbb{X} = \{ \mathbb{X}_{t_\ell}\}_{\ell=1}^{L} = \{ \{\bm{x}_{t_\ell}^k \}_{k=1}^{N_{t_\ell}} \}_{\ell=1}^{L}$, where $\bm{x}_{t_\ell}^k \in \mathbb{R}^d$ is the $k$-th $d$-dimensional sample of state. The aggregate data is actually a $L\times N \times d$ tensor. We refer to $\mathbb{X}_{t_\ell}$ as a snapshot for each $t_\ell$. Our objective is to infer the unknown terms $\bm{m}(\bm{x}), G_{i,j}(\bm{x}), \xi_i(\bm{x}) (i, j = 1, \cdots, d)$ in the stochastic dynamic system directly from the observed aggregate data.

\begin{figure}[tp]
\label{fig:process_chart}
\includegraphics[width=\textwidth]{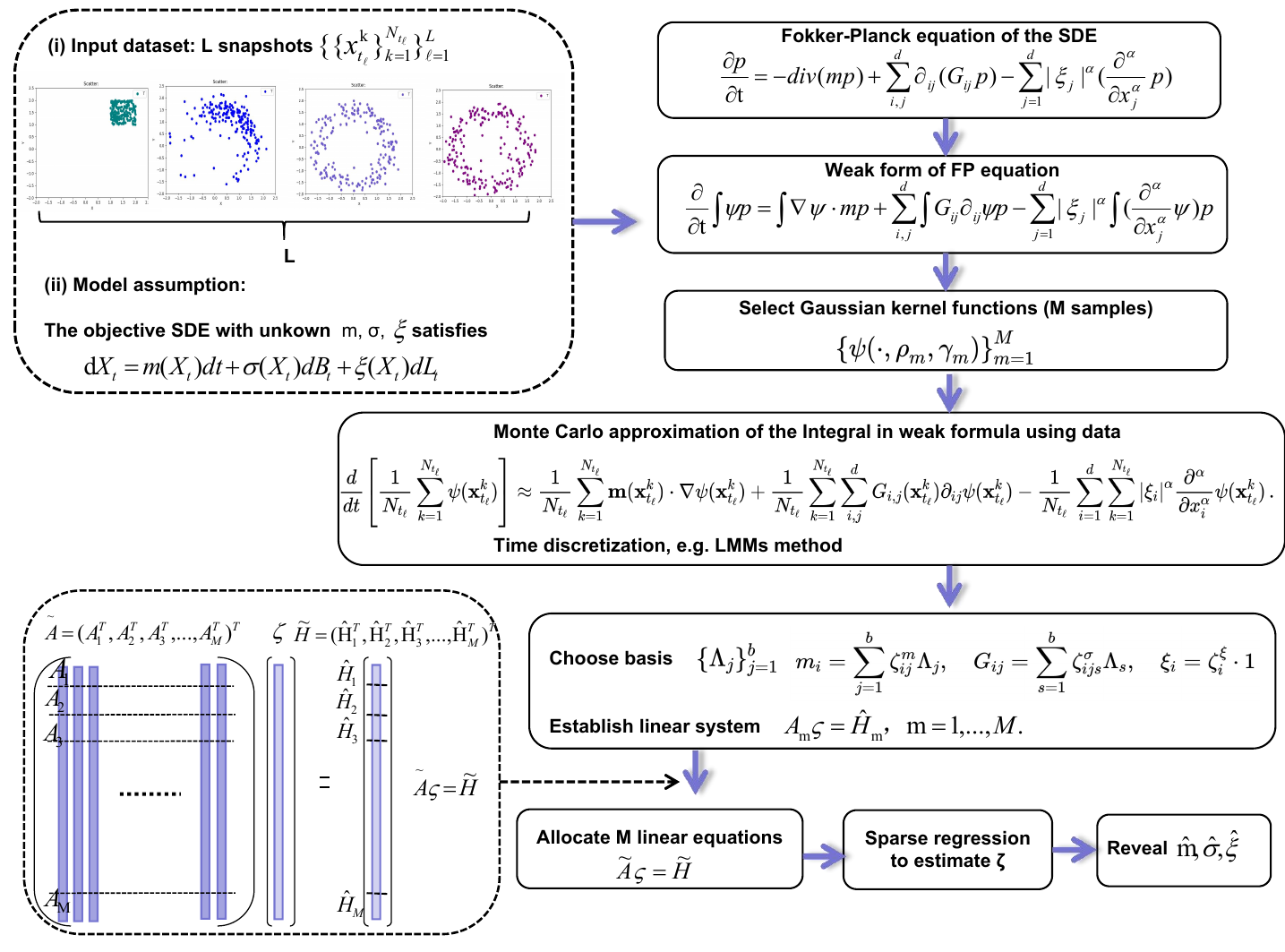}
\caption{Framework of our method. Take the case of independent dimensions as an example. We start by discretely observing the particle motion in a stochastic dynamic system, capturing $L$ snapshots and using the observed particle positions as data set $\mathbb{X}$ (part (i)). Assume that the stochastic system satisfies a SDE with L\'{e}vy noise (part (ii)).
Specifically, the WCR involves the following five steps:
(a) Select Gaussian kernel functions $\psi(\bm{x}, \rho_m, \gamma)$ with the same variance $\gamma$ and utilize the weak form of the Fokker-Planck equation to act the derivatives to the known kernel functions.
(b) Approximate the weak formula using aggregate data by the Monte Carlo method and temporal difference.
(c) Select basis and expand coefficient functions as combinations of basis, then obtain a linear equation system;
(d) Construct a linear equation block $\bm{A}_m \bm{\zeta} = \bm{\hat{H}}_m$ for each sampled Gaussian kernel. They are assembled to a linear equation system of multiple blocks $\bm{\tilde{A}} \bm{\zeta} = \bm{\tilde{H}}$. Solve it by sparse regression;
(e) Estimate the coefficients of the basis expansion by regression. }

\end{figure}

The fact that the data probability density function from a SDE with Gaussian and L\'{e}vy noise is a solution to a Fokker-Planck equation, and they share parameters $\bm{m}(\bm{x})$, $\bm{\sigma}(\bm{x})$ and $\bm{\xi}(\bm{x})$, is a general conclusion that transforms solving stochastic problems into solving deterministic problems. Leveraging this relationship, we consider revealing the unknowns of the FP equation. 
We take the independent case with constant noise intensity $\xi_i$ for L\'{e}vy noise (equation \eqref{equa: fp_ind}) as an example to illustrate the methodology. The framework of our method is summarized as Figure \ref{fig:process_chart} and the main Algorithm \ref{alg:main}.

\begin{algorithm}[htp] 

\renewcommand{\algorithmicrequire}{\textbf{Input:}}
\renewcommand\algorithmicensure {\textbf{Output:} }
\caption{Weak Collocation Regression (WCR) for SDE with L\'{e}vy noise.}

\label{alg:main}

\begin{algorithmic}[1]

\REQUIRE ~~\\

Aggregate data set $\mathbb{X}$ in $d$ dimension, the basis set $\{\Lambda_s\}_{s=1}^b$, $M$ Gaussian kernel functions $\left\{\psi(\cdot, \bm{\rho}_m, \bm{\gamma}) \right\}_{m=1}^M$ ;\\
\ENSURE ~~\\
Estimation of $\bm{m}(\bm{x}), \bm{G}(\bm{x}), \bm{\xi}(\bm{x})$ in equation \eqref{equa:sde_mf_bmlevy};\\

\STATE 
Sample $M$ Gaussian mean values $\{\bm{\rho}_m\}_{m=1}^M$ by LHS method to cover data points $\mathbb{X}$;\\


\FOR{$ m\leq M$}
\STATE Construct linear system of equations $\bm{A}_m \bm{\zeta} = \bm{\hat{H}}_m$ by Algorithm \ref{alg:build_linear}; \\

\ENDFOR

\STATE Stack the linear system $\tilde{\bm{A}}\bm{\zeta} = \bm{\tilde{y}}$
with $\tilde{\bm{A}}=(\bm{A}_1^T,\ldots, \bm{A}_M^T)^T$ and $ \bm{\tilde{H}} = (\bm{\hat{H}}_1^T, \ldots, \bm{\hat{H}}_M^T)^T$;\\

\STATE Compute $\bm{\zeta}$ by ridge regression;\\
\STATE Reconstruct the drift term $\bm{\hat{m}}$, diffusion terms $\bm{\hat{G}}$ and $\bm{\hat{\xi}}$ by expansion with $\bm{\zeta}$. 
\end{algorithmic}
\end{algorithm}

\subsection{Approximate the Weak Form of FP equation by Monte Carlo Method}

\paragraph{Weak formula.}
In order to address the challenge posed by the unattainable derivative of the unknown density function $p(\bm{x}, t)$ in equation \eqref{equa: fp_ind}, we consider the weak form of the Fokker-Planck equation. The existence of the weak formula of FP equation has been established in previous studies \cite{wei2015well,bowles2015weak}, with a simple proof provided in Appendix \ref{appen: weak}. Notably, the weak solution is contingent upon the use of kernel functions, making their selection a crucial consideration. These functions must facilitate the computation of derivatives of various orders, including fractional derivatives, and effectively represent local features such as neighborhoods or distances. Therefore, we opt for Gaussian functions as the kernel functions, which are defined as $\psi(\bm{x}, \bm{\rho}, \bm{\gamma}) \overset{\Delta}{=} \Pi_{i=1}^d \frac{1}{\gamma \sqrt{2 \pi}} e^{-\frac{1}{2}\left(\frac{x_i - \rho_i}{\gamma}\right)^{2}}$, where $\bm{\rho} = \{\rho_i\}_{i=1}^d$ serves as the center of the Gaussian functions, and $\bm{\gamma} = \gamma\cdot \bm{1} \in \mathbb{R}^d$ with $\gamma \in \mathbb{R}$ is the uniform width parameter that controls the radial action range of functions. Furthermore, it is important to note that the weak solution also exhibits smoothness and robustness.
The weak formula for independent $\alpha$-stable L\'{e}vy noises in equation \eqref{equa: fp_ind} is obtained by multiplying kernel function $\psi(\bm{x})$ on both of two sides, integrating on $\bm{x}$, and integration by parts, which is shown in equation \eqref{equa: main_weak}.
\begin{equation}
\begin{aligned}
\label{equa: main_weak}
\frac{d}{d t} \int_{\mathbb{R}^d} p(\bm{x}, t) \psi(\bm{x}) d\bm{x}
    &= \int_{\mathbb{R}^d} \left[m(\bm{x})\cdot \nabla \psi(\bm{x}) \right] p(\bm{x}, t) d\bm{x} + \sum_{i,j}^{d} \int_{\mathbb{R}^d} G_{i,j}(\bm{x}) \partial_{ij}\psi(\bm{x}) p(\bm{x}, t) d\bm{x}\\
    & - \sum_{i=1}^d |\xi_i|^{\alpha}\int_{\mathbb{R}^d} \left(\frac{\partial^\alpha}{\partial x_i^{\alpha}}\psi(\bm{x})\right) p(\bm{x}, t) d\bm{x}\,.
\end{aligned}
\end{equation}

Before dealing with the weak formula, we introduce some useful properties of the fractional Laplacian operator $(-\Delta)^{\alpha/2}$ \cite{burkardt2021unified}, and $-(-\Delta)^{\frac{\alpha}{2}} = \frac{\partial^\alpha}{\partial x^\alpha}$ in one dimension.

\begin{lemma}
\label{lemma: lap_gaussian}
The fractional Laplacian operator has translation and scaling properties. For Schwartz functions $\psi(\bm{x})$, $\bm{x} \in \mathbb{R}^d$ with rapid descent properties, let $\Psi(x) = (-\Delta)^{\alpha/2} \psi(\bm{x})$, and for $\alpha > 0$, the followings hold.
\begin{itemize}
    \item[(i)] $(-\Delta)^{\alpha/2} [\psi(\bm{x} - \bm{x}_0)] = \Psi(\bm{x} - \bm{x}_0)$, for $\bm{x}_0 \in \mathbb{R}^d$;
    \item[(ii)] $(-\Delta)^{\alpha/2} [\psi(\kappa\bm{x})] = \abs{\kappa}^{\alpha}\Psi(\kappa\bm{x})$, for $\kappa \in \mathbb{R}$;
\end{itemize}
Specially, if $\psi(\bm{x})$ is in the form of $e^{-\abs{\bm{x}}^2}$, then for $\alpha \in(0, 1) \cup (1, 2)$, the fractional Laplacian of $\psi(\bm{x})$ can be calculated by
\begin{equation}
    \begin{aligned}
    (-\Delta)^{\alpha/2} \psi(\bm{x}) = \frac{2^{\alpha} \Gamma((d+\alpha)/2)}{\Gamma(d/2)} \,
    _1F_1(\frac{d+\alpha}{2}; \frac{d}{2}; -\abs{\bm{x}}^2)\,, \quad \bm{x} \in \mathbb{R}^d\,,
    \end{aligned}
\end{equation}
where $_1F_1$ is the confluent hypergeometric function. 
\end{lemma}

Therefore, the fractional derivative of the $i$-th component of Gaussian kernel $\psi(\bm{x}) = \psi(\bm{x}, \bm{\rho}, \bm{\gamma})$ can be solved by the $_1F_1$ function.
\begin{equation}
    \begin{aligned}
\frac{\partial^\alpha}{\partial x_i^{\alpha}}\psi(\bm{x}) 
    =& \frac{1}{(\sqrt{2\pi}\gamma)^d}\exp\left[-\frac{1}{2\gamma^2} \sum \limits_{j=1, j\neq i}^d (x_j - \rho_j)^2 \right] \frac{2^{\alpha} \Gamma((1+\alpha)/2)}{\Gamma(1/2)} \abs{\frac{1}{\sqrt{2} \gamma}}^{\alpha}\, _1F_1 \left(\frac{1+\alpha}{2}; \frac{1}{2}; -\frac{ (x_i - \rho_i)^2}{2\gamma^2}\right)\,.
    \end{aligned}
\end{equation}

Lemma \ref{lemma: lap_gaussian} provides a solution for the fractional derivative, then only unknown terms in equation \eqref{equa: main_weak} are $\bm{m}$, $G_{i,j}$, $\xi_i$, and the density function $p(\bm{x}, t)$. Recall that our work does not require obtaining data distribution; we aim to solve unknown parameters and obtain the modeled SDE explicitly. We can approximate the integration of the probability density function $p(\bm{x}, t)$ using the aggregate data through Monte Carlo (MC) method.

\paragraph{Temporal derivative and Monte Carlo approximation.}
We begin by rewriting equation \eqref{equa: main_weak} due to the average property of the probability density function:
\begin{equation}
    \begin{aligned}
        \frac{d}{dt} \mathbb{E}_{\bm{x} \sim p(\bm{x},t)} [\psi(\bm{x})] = \mathbb{E}_{\bm{x} \sim p(\bm{x},t)} \left[ \bm{m}(\bm{x})\cdot \nabla \psi(\bm{x}) + \sum_{i,j}^{d} G_{i,j}(\bm{x}) \partial_{ij}\psi(\bm{x}) - \sum_{i=1}^d |\xi_i|^\alpha \frac{\partial^\alpha}{\partial x_i^\alpha} \psi(\bm{x}) \right]\,.
    \end{aligned}
\end{equation}
where $\mathbb{E}_{\bm{x}\sim p(\bm{x},t)}$ is the expectation over the probability density function. To approximate the expectation terms, we use the Monte Carlo method, known for its robustness and accuracy in sampling, with a dimension-independent convergence rate \cite{Lemieux2009MonteCA}. The advantage of employing the Monte Carlo method is its ability to avoid complex integration calculations. It is worth noting that at time $t_{\ell}$, the $k$-th sample is represented as $x_{t_{\ell}}^k$, and we collect $N_{t_{\ell}}$ samples at this time. Using the relationship $\mathbb{E}_{\bm{x}\sim p(\bm{x}, t)}[\psi(\bm{x})] \approx \frac{1}{N_{t_\ell}} \sum\limits_{k=1}^{N_{t_\ell}} \psi(\bm{x}_{t_\ell}^{k})$ for $\psi$, we can rewrite equation \eqref{equa: main_weak}.
\begin{equation}
\begin{aligned}
\label{equa: monte_carlo}
\frac{d}{dt} \left[\frac{1}{N_{t_\ell}} \sum_{k=1}^{N_{t_\ell}} \psi(\bm{x}_{t_\ell}^{k}) \right] 
& \approx \frac{1}{N_{t_\ell}} \sum_{k=1}^{N_{t_\ell}} \bm{m}(\bm{x}_{t_\ell}^{k}) \cdot \nabla\psi(\bm{x}_{t_\ell}^{k}) + \frac{1}{N_{t_\ell}} \sum_{k=1}^{N_{t_\ell}} \sum_{i,j}^{d} G_{i,j} (\bm{x}_{t_\ell}^{k}) \partial_{ij} \psi(\bm{x}_{t_\ell}^{k}) - \frac{1}{N_{t_\ell}} \sum_{i=1}^d \sum_{k=1}^{N_{t_\ell}}|\xi_i|^{\alpha} \frac{\partial^\alpha}{\partial x_i^{\alpha}}\psi(\bm{x}_{t_\ell}^{k}) \,.
\end{aligned}
\end{equation}

The left side of the above formula still contains a temporal derivative, and we approximate it by the linear multi-step methods (LMMs) \cite{du2022discovery}. Recall that $\bm{x}_{t_\ell} = (x_{t_\ell}^1, \cdots, x_{t_\ell}^N)$, and we denote equation \eqref{equa: monte_carlo} as $\frac{d}{dt} h(\bm{x}_{t_\ell}) = s(\bm{x}_{t_\ell})$, where $h(\bm{x}_{t_\ell})$ and $s(\bm{x}_{t_\ell})$ are scalar functions.
We take trapezoidal formula for time discretization as an example.
\begin{equation}
    \begin{aligned}
    \label{equa: LMM}
        h(\bm{x}_{t_{\ell+1}}) - h(\bm{x}_\ell) &= \frac{\Delta t}{2} \left( s(\bm{x}_{t_\ell})+ s (\bm{x}_{t_{\ell+1}}) \right) \,, \quad \ell = 1,\cdots, L-1 \,.
    \end{aligned}
\end{equation}

\subsection{Build Linear Regression Equations}
\label{sec: build_linear}
This section focuses on selecting basis and constructing regression equations to solve the discrete form of equation \eqref{equa: monte_carlo} using the "collocation method". This method simplifies the solving of differential equations by transforming them into a system of linear equations. Specifically, we consider sampling data at irregular time and space grids.

\paragraph{Selected basis.}
We choose basis $\{\Lambda_s\}_{s=1}^b$ for the drift term and the diffusion term of Gaussian noise based on states $\bm{x}$, represented as
\begin{equation}
\label{equa: basis_wcr_1}
\Lambda = \{\Lambda_1, \Lambda_2, \cdots, \Lambda_b \}^T \,,
\end{equation}
where $b$ is the number of the elements in $\Lambda$. Simultaneously, the coefficient matrix $\bm{\xi}$ is a diagonal matrix with constant elements $\xi_i, i=1, 2, \cdots, d$, with the basis $\{1, 1,\cdots, 1\}_d$. Notably, the parameters of each dimension $\xi_i$ may differ, necessitating learning $d$ parameters. Subsequently, the components $m_i$ of $\bm{m}(\bm{x})=[m_i]$, $G_{ij}$ of $\bm{G}(\bm{x})=[G_{ij}]$, and $\xi_i$ of $[\xi_i]$ are expanded as
\begin{equation}
\label{eq.expan}
m_i = \sum_{j=1}^{b} \zeta_{ij}^m \Lambda_j, \quad G_{ij} = \sum_{s=1}^{b} \zeta_{ijs}^{\sigma} \Lambda_s, \quad \xi_i = \zeta_i^{\xi} \cdot 1\,,
\end{equation}
with $\Lambda_s$ representing the $s$-th component in the basis set $\Lambda$, and $\zeta_{ij}^m$, $\zeta_{ijs}^\sigma$, $\zeta_i^\xi$ acting as coefficients of the expansion of these terms. 
Specifically, for $d$-dimensional variable $\bm{x}$, we have
\begin{equation}
\begin{aligned}
\label{equa:expand_m}
\bm{m}(\bm{x}) =
\begin{pmatrix}
m_1(\bm{x})\\
m_2(\bm{x})\\
\vdots\\
m_d(\bm{x})
\end{pmatrix}
=
\begin{pmatrix}
\sum_{j=1}^{b} \zeta_{1j}^m \Lambda_j\\
\sum_{j=1}^{b} \zeta_{2j}^m \Lambda_j\\
\vdots\\
\sum_{j=1}^{b} \zeta_{dj}^m \Lambda_j
\end{pmatrix}
=
\begin{pmatrix}
\zeta_{11}^m, &\cdots, &\zeta_{1b}^m \\
\zeta_{21}^m, &\cdots, &\zeta_{2b}^m \\
&\ddots & \\
\zeta_{d1}^m, &\cdots, &\zeta_{db}^m
\end{pmatrix} \cdot \begin{pmatrix}
\Lambda_1\\
\Lambda_2\\
\vdots\\
\Lambda_b
\end{pmatrix}\,.
\end{aligned}
\end{equation}

The coefficient matrix of the drift term $\{\zeta_{ij}^m\}_{db}$, takes the form of a matrix. The coefficient matrix of $\bm{G}$, $\{ \zeta_{ijs}^{\sigma}\}_{d\times d\times b}$ is actually a tensor in $\mathbb{R}^{d\times d\times b}$, as each element $G_{i,j}(\bm{x})$ of $\bm{G}$ can be represented by a linear combination of a $b \times 1$ dimensional basis. We flatten the matrix or tensor into a vector according to the dimension $d$, so that equation \eqref{equa: monte_carlo} is expressed as the inner product of the basis expansion coefficient and the basis. Let
\begin{equation} 
\begin{aligned}
\label{equa: model_basis}
    \bm{\zeta} = &\{ \underbrace{\zeta_{11}^m, \cdots, \zeta_{1b}^m, \cdots, \zeta_{d1}^m, \cdots, \zeta_{db}^m}_{\text{flatten of }\bm{\zeta}^m}, \, \underbrace{\zeta_{111}^{\sigma}, \cdots, \zeta_{11b}^{\sigma}, \cdots, \zeta_{dd1}^{\sigma}, \cdots, \zeta_{ddb}^{\sigma}}_{\text{flatten of }\bm{\zeta}^{\sigma}}, \, \underbrace{\zeta_{1}^{\xi}, \cdots, \zeta_{d}^{\xi}}_{\text{flatten of }\bm{\zeta}^{\xi}}\}^T \in \mathbb{R}^{d^2b+db+d}\,,
    \end{aligned}
\end{equation}
be the flatten coefficient vector of the combination of the three basis set, then equation \eqref{equa: monte_carlo} is transformed into 
\begin{equation}
\begin{aligned}
\label{equa: expand_basis}
\frac{d}{dt} \left[\frac{1}{N_t} \sum_{k=1}^{N_t} \psi(\bm{x}_{t}^{k}) \right]
 &= \sum_{j=1}^{b} \frac{1}{N_t} \sum_{k=1}^{N_t} 
   \sum_{i=1}^{d} \psi_{x_i} (\bm{x}_{t}^{k}) \Lambda_j \zeta_{ij}^m + \sum_{s=1}^{b} \frac{1}{N_t} \sum_{k=1}^{N_t} \sum_{i,j}^{d} \partial_{ij} \psi(\bm{x}_{t}^{k}) \Lambda_s \zeta_{ijs}^{\sigma} - \frac{1}{N_t} \sum_{k=1}^{N_t} \sum_{i}^{d} \frac{\partial^\alpha}{\partial x_i^{\alpha}} \psi(\bm{x}_{t}^{k})\zeta_{i}^{\xi}\\
   &\triangleq \bm{B}(\bm{x}_t) \cdot \bm{\zeta}\,.
\end{aligned}
\end{equation}

We begin by discretizing the approximation equation \eqref{equa: model_basis} at observation time $\ell \in \{1, 2, \cdots, L\}$ using the LMMs in equation \eqref{equa: LMM}. For each $\ell\in \{1, \cdots, L-1\}$, the temporal discrete format is given by 
\begin{equation}
    \begin{aligned}
    \label{equa: LMM2}
        \hat{h}_\ell (\mathbb{X}) &= \hat{h}(\bm{x}_{t_\ell}) = h(\bm{x}_{t_{\ell+1}}) - h(\bm{x}_{t_\ell}) = \left(\frac{1}{N_{t_{\ell+1}}} \sum_{k=1}^{N_{t_{\ell+1}}} \psi(\bm{x}_{t_{\ell+1}}^{k})\right) - \left( \frac{1}{N_{t_{\ell}} }\sum_{k=1}^{N_{t_{\ell}} }\psi(\bm{x}_{t_\ell}^{k})\right) \\
&=\frac{\Delta t}{2} (\bm{b}_\ell^T(\bm{x}_{t_\ell}) +\bm{b}_{\ell+1}^T(\bm{x}_{t_{\ell+1}})) \cdot \bm{\zeta} \,,
    \end{aligned}
\end{equation}
where $\bm{b}_\ell^T$ is the $\ell$-th row of $\bm{B}$ in equation \eqref{equa: model_basis}, and $\bm{x}_{t_\ell}$ is the $\ell$-th observation in dataset $\mathbb{X}$. 
Next, we assemble these $L$ equations. Let $\bm{\hat{H}}(\mathbb{X}) = \left(\hat{h}_1(\mathbb{X}), \cdots, \hat{h}_L(\mathbb{X}) \right)^T \in \mathbb{R}^L$, and $\bm{a}_\ell^T(\mathbb{X}) = \frac{\Delta t}{2} (\bm{b}_\ell^T(\bm{x}_{t_\ell}) +\bm{b}_{\ell+1}^T(\bm{x}_{t_{\ell+1}}))$, $\ell = 1, \cdots, L-1$. Consequently, for$ \bm{A}(\mathbb{X}):= \begin{pmatrix}
\bf{a}^T_1(\mathbb{X})\\ \vdots\\\bf{a}^T_{L-1}(\mathbb{X})\end{pmatrix} \in \mathbb{R}^{L\times(db + d^2b+ d)}$, we obtain $\bm{A}(\mathbb{X}) \bm{\zeta} = \bm{\hat{H}}(\mathbb{X})$, which gives the relationship between the data set $\mathbb{X}$ and the parameter of basis expansions for the hidden dynamics. We can see that WCR transforms the problem of solving PDE into a system of linear equations. For pseudo-code, see Algorithm \ref{alg:build_linear}.

\begin{algorithm}[htb] 

\renewcommand{\algorithmicrequire}{\textbf{Input:}}

\renewcommand\algorithmicensure {\textbf{Output:} }

\caption{Build Linear System (for each $\bm{A}_m$, $m=1,2,\cdots,M$)}
\label{alg:build_linear}

\begin{algorithmic}[1]

\REQUIRE ~~\\
Aggregate dataset $\mathbb{X}$, the kernel function $\psi(\bm{x}, \bm{\rho}, \gamma)$, and the basis set $\{\Lambda_s\}_{s=1}^b$;
\ENSURE ~~\\
A linear equation system $\bm{A}\bm{\zeta} = \bm{\hat{H}}$;\\

\STATE Monte Carlo approximation of weak formula (eq.\eqref{equa: main_weak}) by data points;
\STATE Obtain derivative information of the kernel function;\\
\STATE Expand the drift vector $\bm{m}=[m_i]$, diffusion matrix $\bm{G}=[G_{ij}]$, and $\bm{\xi} = diag\{ \xi_1, \cdots, \xi_d\}$ by basis combination (eq.\eqref{equa: basis_wcr_1})
$
m_i = \sum_{j=1}^{b} \zeta_{ij}^m \Lambda_j, G_{ij} = \sum_{k=1}^{b} \zeta_{ijk}^{\sigma} \Lambda_k,
\text{and } \xi_i = \zeta_{i}^{\xi}
;$ \\
\STATE Time discretization by LMM (eq.\eqref{equa: LMM2});
\STATE Assemble $L$ time discretization equations and obtain a linear equation system (eq.\eqref{equa:linear});\\

\end{algorithmic}
\end{algorithm}

\paragraph{Assembly of block linear equations.}
The linear system for one selected Gaussian kernel $\psi(\bm{x}, \bm{\rho}, \bm{\gamma})$ is
\begin{equation}
\label{equa:linear}
\bm{A}(\mathbb{X}, \bm{\rho}, \bm{\gamma}) \bm{\zeta} = \bm{\hat{H}}(\mathbb{X}, \bm{\rho}, \bm{\gamma}) \,.
\end{equation}
However, one type of Gaussian function has low efficiency that it can only learn information from the neighborhood of itself. To cover all data points more fully, we sample $M$ Gaussian functions with distinct mean values $\rho_m, m=1, \cdots, M$ based on data distribution, and repeat building linear equation as equation \eqref{equa:linear} for $M$ times.
The mean value $\rho_m$ is sampled by Latin Hypercub Sampling (LHS) method, which is one of Monte Carlo-based approaches and reduces the variance of an estimator compared to random sampling \cite{iman2008atin}. It leads to a reduction of the sample size while maintaining the statistical significance.
Similar to equation \eqref{equa:linear}, the $m$-th Gaussian function $\psi(\cdot, \bm{\rho}_m, \bm{\gamma})$ yields a linear equation system $\bm{A}_m\bm{\zeta} =\bm{\hat{H}}_m$, where $A_m = \bm{A}(\mathbb{X}, \bm{\rho}_m, \bm{\gamma})$ and $\bm{\hat{H}}_m  =  \bm{\hat{H}}(\mathbb{X}, \bm{\rho}_m, \bm{\gamma})$. We build the linear system over different Gaussian functions as equation \eqref{equa: large} by assembling $M$ blocks, intuitively represented by Figure \ref{fig:process_chart}.
\begin{equation}
\label{equa: large}
\bm{\tilde{A}} \bm{\zeta} = \bm{\tilde{H}},
\end{equation}
where $\bm{\tilde{A}} = \begin{pmatrix}
 \bm{A}_1^T, \bm{A}_2^T, \cdots, \bm{A}_M^T 
\end{pmatrix}^T \in \mathbb{R}^{ML \times (db+d^2b+d)},\, \text{and} \,
\bm{\tilde{H}} = 
\begin{pmatrix}
 \bm{\hat{H}}_1^T, 
 \bm{\hat{H}}_2^T, 
 \cdots,
 \bm{\hat{H}}_M^T
\end{pmatrix}^T \in \mathbb{R}^{ML\times 1}$.

\paragraph{Solution of unknown parameters.}
From equation \eqref{equa: large}, we obtain the estimation of the unknown vector $\bm{\zeta}$ by applying sequential threshold ridge regression (STRidge) technique \cite{brunton2016discovering}. 
The data we collect often exhibits characteristics of sparsity. To address this, we incorporate a $2$-norm regularization term within the framework of regression. For the system of linear equations, we employ least square regression to obtain the estimated $\hat{\bm{\zeta}}$. We further optimize the obtained solution. The STRidge algorithm is initiated with the least squares solution and systematically sets all coefficients falling below a predetermined threshold to zero. This operation is repeated iteratively, with the indices of non-zero coefficients identified in each iteration. The process continues until the non-zero coefficients converge, effectively determining the optimal solution that best represents the underlying dynamics.
By distributing the components of $\hat{\bm{\zeta}}$ to drift and diffusion terms, the unknown term $\bm{m}(\bm{x}), \bm{G}(\bm{x}), \bm{\xi}$ are obtained by combinations of $\hat{\bm{\zeta}}$ and basis. The whole process can be simplified to Algorithm \ref{alg:main}.

\begin{remark}
\label{remark:RI}
When the L\'{e}vy process is not independent in each dimension, the L\'{e}vy measure is not equal to the sum of the measure of each dimension.
Specially, when the $\alpha$-stable L\'{e}vy process ($\alpha \in (0, 1) \cup (1, 2)$) is \textit{rotation invariant} (RI)\footnote{A stochastic process $L_t$ is \textit{rotation invariant} if $L_t \overset{\text{d}}{=} U L_t$ for every orthogonal matrix $U$, where $\overset{\text{d}}{=}$ means equally distributed on both sides of the equation.}, the Fokker-Planck of the SDE $d\bm{X}_t = \bm{m}(\bm{X}_t) dt + \bm{\sigma}(\bm{X}_t) d \bm{B}_t + \bm{\xi}_0(\bm{X}_t) d \bm{L}_t$ is \cite{duan2015introduction}
\begin{equation}
    \begin{aligned}
    \label{equa: fp_rotation}
    \frac{\partial}{\partial t} p(\bm{x}, t) 
    &= -div  \left(\bm{m}(\bm{x}) p(\bm{x}, t)\right) + \sum_{i,j}^{d} \partial_{ij} \left( G_{i,j}(\bm{x}) p(\bm{x}, t) \right) - \theta_{\alpha, d} (-\Delta)^{\alpha/2} \left(|\xi_0(\bm{x})|^{\alpha} p(\bm{x}, t)\right) \,
\end{aligned}
\end{equation}
where $\theta_{\alpha, d} = \pi^{-1/2} \frac{\Gamma(\frac{1+\alpha}{2}) \Gamma(\frac{d}{2})}{\Gamma(\frac{d+\alpha}{2})}$, and $\bm{\xi}_0(\bm{x}) = \xi_0(\bm{x})\cdot \bm{I}$ with $\xi_0(\bm{x})$ is a scalar function and $\bm{I}$ is an unit matrix.
\end{remark}

\section{Experimental Results}
\label{sec: exper}

\subsection{Experimental Setting}

\paragraph{Data.} 
The aggregate data set, denoted as $\mathbb{X} = \{ \mathbb{X}_{t_\ell}\}_{\ell=1}^{L} = \{ \{\bm{x}_{t_\ell}^k \}_{k=1}^{N_{t_\ell}} \}_{\ell=1}^{L}$, is composed of snapshots at discrete times $\ell = 1, \cdots, L$. These snapshots capture data positions $ \{\bm{x}_{t_\ell}^k\}_{k=1}^{N_{t_\ell}}$, where $N_{t_\ell}$ signifies the number of data points at time $t_\ell$. The observed data can be experimental data, and in our experiments, they are generated snapshots based on model assumptions with an initial data distribution of $\cN(0, 0.2)$, where $\cN$ represents the normal distribution. Details of the data generation process are provided in Appendix \ref{appen: data}, with data being generated within the time range of $t=0, 0.1, \cdots, 1.8$ seconds. All training data used in the experiments are from the time range of $t \in [0, 1.0]$, while data after $t=1.0$ are used to test the performance.

\paragraph{Basis and kernel functions.} The process of converting the weak form of a FP equation into a system of linear equations is elucidated in the preceding section. This is achieved through the construction of a basis set and the selection of kernel functions. Due to their ability to locally approximate a smooth function, as noted by \cite{ganzburg2008polynomial}, polynomials are chosen as the basis function for collocation methods. This choice is driven by their simplicity and ease of operation. The basis functions for the $p$-th order polynomial with respect to the variable $\bm{x}=(x_1,\cdots, x_d)^T$ are expressed as follows:
\begin{equation}
\label{equa: basis_wcr}
\Lambda = \{1, x_1, x_2, \cdots, x_d, x_1x_1, x_1x_2, \cdots, x_1 x_p, x_2 x_2, \cdots, x_d^p\}^T \,.
\end{equation}
The number of the terms of $\Lambda$ is $b \overset{\Delta}{=}|\Lambda| = \begin{pmatrix}
p+d\\
p
\end{pmatrix}$.
We choose $M$ Gaussian functions as kernel functions all possessing the same variance $\bm{\gamma} = \gamma \bm{I}_d$ and mean values $\rho_m$, $m = 1, 2, \cdots, M$. These mean values, which are contingent on the data distribution, are sampled between the minimum value $l_b = \min\limits_{\ell=1, \cdots, L}\min\limits_{k=1, \cdots, N_{t_\ell}} \bm{x}_{t_\ell}^{k}$ and maximum value $u_b = \max\limits_{\ell=1, \cdots, L}\max\limits_{k=1, \cdots, N_{t_\ell}} \bm{x}_{t_\ell}^{k}$ of the data points to encompass all data points. The mean values are sampled by Latin Hypercube Sampling (LHS) method \cite{McKay1992lhs}. After choosing basis and building a linear system $\bm{\tilde{A}} \bm{\zeta} = \bm{\tilde{H}}$, we apply the sequential threshold ridge regression (STRidge) technique \cite{brunton2016discovering} to obtain the estimated $\hat{\bm{\zeta}}$. 

\paragraph{Error.} To measure the error, we use the Maximum Relative Error (MRE) of non-zero terms, which is defined as 
\begin{equation*}
    \text{MRE} = \max_{\theta_i\neq 0}\frac{|\hat{\theta}_i-\theta_i|}{|\theta_i|}\,, 
\end{equation*}
where $\theta_i$ represents the $i$-th parameter of the drift and diffusion terms and $\hat{\theta}_i$ is the estimated parameter from data. 

We illustrate the importance of considering L\'{e}vy noise in Section \ref{sec: necessity}, and subsequent sections serve to substantiate the efficacy of our approach in various situations, encompassing comparison with the previous work (Section \ref{subsec:1d_ind}), multi-dimensional independent L\'{e}vy noise across dimensions, coupled cases (Section \ref{sec:multi} and \ref{sec:5d}), non-polynomial drift term configurations (Section \ref{sec: nonpoly}), functional diffusion examples (Section \ref{sec:func_diff}). Lastly, Section \ref{sec:robust} tests the robustness. All experiments are implemented on laptop with 11th Gen Intel(R) Core(TM) i7-1165G7 @ 2.80GHz. Our code will be open-sourced later.

\subsection{Necessity of SDE Model with L\'{e}vy Noise}
\label{sec: necessity}

\begin{figure}[t]
\centering    
\includegraphics[width=5.2cm,height=3.7cm]{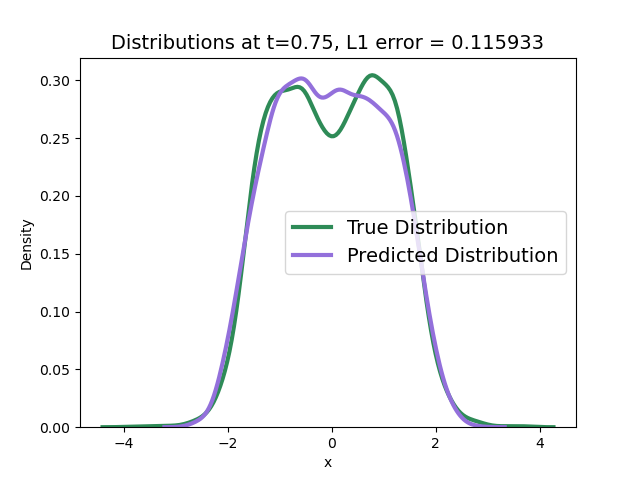}
\includegraphics[width=5.2cm,height=3.7cm]{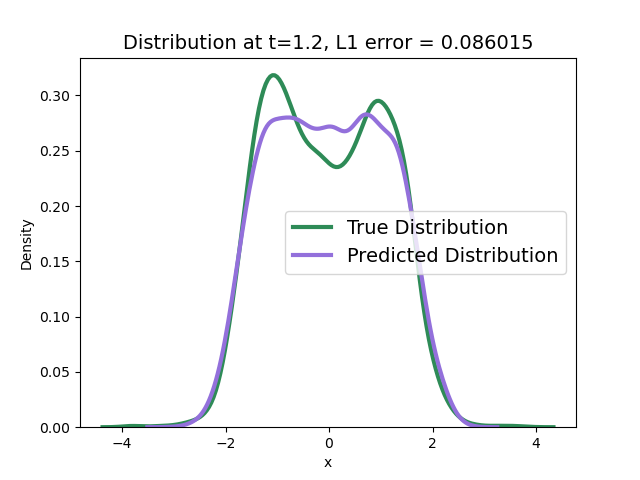}

\includegraphics[width=5.2cm,height=3.7cm]{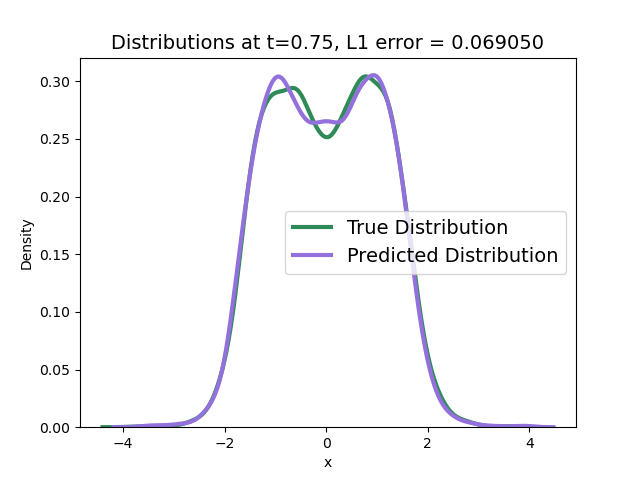}
\includegraphics[width=5.2cm,height=3.7cm]{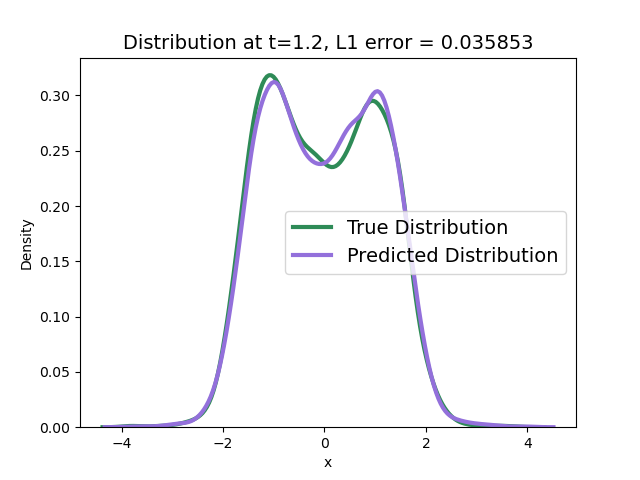}
\caption{Comparison of the distribution of true data and the predicted distributions of estimated SDE with or without L\'{e}vy noise.
The top row represents results obtained from WCR only with Gaussian noise and the bottom row shows the predicted distribution by WCR method based on both Gaussian and L\'{e}vy noise.
The data distributions are plotted at $t=0.75$ (left) and $t=1.2$ (right).}
\label{fig: why_levy}
\end{figure}

Stochastic systems often exhibit heavy-tailed properties or jumps in financial decision-making or weather forecasting \cite{emmer2004optimal, zheng2020maximum}, the WCR based solely on Gaussian noise becomes inadequate for accurately capturing the data evolution. We model a SDE with L\'{e}vy noise. Initially, we utilize WCR only with Gaussian noise in \cite{lu2024weak}, in the form of $d\bm{X}_t = \bm{m}(\bm{X}_t) dt + \bm{\sigma}(\bm{X}_t) d\bm{B}_t$, to estimate the more complex system. The estimated SDE is obtained by substituting the unknown parameters, obtained by WCR from the training data $t=\{0, 0.1, \cdots, 1.0\}$, into the modeled SDE. To evaluate the accuracy of the estimates, we choose $t=0.75$ and $t = 1.2$ as the test moments. At $t = 1.2$, we compare the true distribution of the discrete test dataset with the predicted distribution of the estimated SDE evolution to the moment $t=1.2$. Additionally, at $t = 0.75$, we compare the interpolated results of the real data with the predicted distribution of the estimated SDE evolved to $t=0.75$. 
The first row of Figure \ref{fig: why_levy} displays the data distributions at $t=0.75$ and at $t=1.2$ using density profiles of the data tensor at the test moments fitted with a $kde$ in python. However, the results only considering Gaussian noise are not optimal, as indicated by the labeled $L_1$ error, and they do not recognize bimodal peaks well. 

Additionally, we consider the system with both Gaussian noise and L\'{e}vy noise and apply the previously described WCR method to estimate the SDE. Results (the bottom row of Figure \ref{fig: why_levy}) reveal that accounting for L\'{e}vy noise extends the method's applicability to a broader spectrum of problems. Despite the existence of numerous noise types in real world, the incorporation of L\'{e}vy noise could cover a diverse range of problems. Our experiments sets the stable index $\alpha$ of L\'{e}vy process as $1.5$. In more general scenarios without fixed $\alpha$, the estimation of $\alpha$ may be necessitated \cite{fang2022end}.

\subsection{Comparison to the Existing Work}
\label{subsec:1d_ind}

Considered as a one-dimensional stochastic model, previous methods \cite{yang2020generative,chen2021solving} have demonstrated promising results in learning stochastic differential equations (SDE) that integrate both Gaussian and non-Gaussian noise. But screening out Gaussian noise at relative distant locations is difficult due to the different diffusion rates of L\'{e}vy and Gaussian processes.
We consider a one-dimensional stochastic model as
\begin{equation}
    \begin{aligned}
        \label{equa: exp2}
        d X_t = (X_t - X_t^3) dt + d B_t + d L_t \,,
    \end{aligned}
\end{equation}
The basis for the drift term is chosen as $\Lambda = \{1, x, x^2, x^3\}$, and our aim is to reveal unknown parameters $\bm{\zeta} = (\lambda_0, \lambda_1, \lambda_2, \lambda_3, \sigma, \xi)^T$ in 
$d X_t = (\lambda_0 + \lambda_1 X_t + \lambda_2 X_t^2 + \lambda_3 X_t^3) dt + \sigma dB_t + \xi d L_t \,.$

\begin{table}[thp]
\centering 
\scalebox{0.9}{
\begin{tabular}{ccccccc}
\toprule Parameter & $\lambda_{0}$ & $\lambda_{1}$ & $\lambda_{2}$ & $\lambda_{3}$ & $\sigma$ & $\xi$ \\
\hline
True parameters & 0 & 1 & 0 & $-1$ & 1 & 1 \\
\hline
(a) \cite{chen2021solving} & $0.0790$ & $\bm{0.9860}$ & $-0.0314$ & $-0.8644$ & $0.6966$ & $\bm{1.0318}$ \\
(a) WCR & $\bm{0.0}$ & $0.8585473$ & $\bm{0.0}$ & $\bf{-0.9083656}$ & $\bf{1.0418}$ & $0.8964$ \\
\hline
(b) \cite{chen2021solving} & $0.0042$ & $\bm{0.9391}$ & $-0.0380$ & $\bm{-0.9960}$ & $0.7076$ & $1.2092$ \\
(b) WCR & $\bm{0.0}$ & $0.9140$ & $\bm{0.0}$ & $-0.9585$ & $\bm{0.9221}$ & $\bm{1.0262}$ \\
\hline
(c) WCR with more snapshots & $0.0$ & $0.9796$ & $0.0$ & $-1.0141$ & $1.0167$ & $0.9948$ \\
\bottomrule
\end{tabular}
}
\caption{Comparison with \cite{chen2021solving} of the $1d$ results of revealing the unknown terms with $10,000$ samples at different snapshots. (a) Observations at $t=0.2, 0.5, 1.0$ s; (b) Observations at $t=0.1, 0.4, 0.7, 1.0$ s; (c) Observations at $t=0.1, 0.2, 0.3, \cdots, 0.9, 1.0$ s; All settings have the same initial distribution $\cN(0, 0.2)$.} 
\label{table: both_1d_compare}
\end{table}

In Table \ref{table: both_1d_compare}, the WCR algorithm is compared with a previous physics-informed neural networks (PINNs) based method \cite{chen2021solving} to demonstrate its advantages. Our results show that for the same discrete observation moments ($t = 0.2, 0.5, 1.0$ and $t = 0.1, 0.4, 0.7, 1.0$ respectively) of the same SDE \eqref{equa: exp2}, the WCR algorithm achieves lower MRE under the same experimental conditions, including the initial data distribution (a Gaussian distribution $\cN(0, 0.2)$) and a sample size of $10,000$. This indicates that our approach is not only more time-efficient but also offers improved accuracy. Specifically, our method requires approximately $1.5 s$, $1.0 s$, and $7.8 s$ of computational time on \textit{CPU} device for cases $(a), (b), (c)$, respectively. This time efficiency can be attributed to two primary factors. Firstly, we adopt a simple linear regression as opposed to a more complex neural network. Secondly, the FP equation contains an integral term, which is not only complicated to approximate but also challenging to extend to higher dimensions. We consider Lemma \ref{lemma: lap_gaussian}, thereby reducing the computational burden associated. Furthermore, WCR exhibits reduced error with an increase in the number of snapshots, which is influenced by the temporal difference error.
Their work conducts experiments with L\'{e}vy noise in a two dimensional case, and we test in higher dimensions in later context.

\subsection{Multi-dimensional Problems}
\label{sec:multi}

\subsubsection{Independent Across Dimensions}
\label{sec: ind_multi}

We consider the multi-dimensional cases where the drift terms are cubic polynomials $x_i - x_i^3$, $i=1,\cdots,d$, the diffusion terms are constants, and $\bm{\xi} (\bm{X_t}) = diag\{\xi_1, \cdots, \xi_d\}$ with $\xi_i$ are constants.
\begin{equation}
    \begin{aligned}
    \label{equa:indep_dim2}
        d X_{i,t} &= (\lambda_0^{(i)} + \lambda_1^{(i)} X_{i,t} + \lambda_2^{(i)} X_{i,t}^2 + \lambda_3^{(i)} X_{i,t}^3) dt + \sigma^{(i)} d B_{i,t} + \xi^{(i)} d L_{i,t}
    \end{aligned}
\end{equation}

\begin{example}
\label{exam: 2d_ind}
We first consider independent two dimensional problem with drift terms $m_i(x_i) = x_i - x_i^3$ and noise intensities $\sigma^{(i)}=1, \xi^{(i)}=1$, $i=1,2$. 
We sample $M=150$ Gaussian kernels.
Table \ref{table: 2d_levy_inequal_Diag} shows the estimated results. 
\end{example}

\begin{table}[htp]
\centering
\scalebox{0.9}{
\begin{tabular}{ccccccc}
\toprule Parameter & $\lambda_{0}^{(1)}$ & $\lambda_{1}^{(1)}$ & $\lambda_{2}^{(1)}$ & $\lambda_{3}^{(1)}$ & $\sigma^{(1)}$ & $\xi^{(1)}$ \\
\hline
True & 0 & 1 & 0 & $-1$ & 1 & 1 \\
\hline
Estimated & $-0.08049$ & $0.93940$ & $-0.12519$ & $-1.06913$ & $1.10025$ & $0.90092$ \\
\hline
Parameter & $\lambda_{0}^{(2)}$ & $\lambda_{1}^{(2)}$ & $\lambda_{2}^{(2)}$ & $\lambda_{3}^{(2)}$ & $\sigma^{(2)}$ & $\xi^{(2)}$ \\
\hline
True & 0 & 1 & 0 & $-1$ & 1 & 1 \\
\midrule
Estimated & $0.0$ & $1.10810$ & $0.0$ & $-1.11434$ & $1.09796$ &  $0.985627$ \\
\bottomrule
\end{tabular}
}
\caption{Reveal the unknown drift and diffusion terms with $10,000$ samples of $X_t$ at time snapshots $t = \{0.1, 0.3, 0.5, 0.7, 1.0\}$ in $2d$ case of Example \ref{exam: 2d_ind}.} 
\label{table: 2d_levy_inequal_Diag}
\end{table}
The increase in dimension and sample space results in more Gaussian kernel samplings and causes more computational burdens.

\begin{example}
Our method can also distinguish various noises in a multi-noise system. For example, consider a two-dimensional dynamic system with distinct noises in each dimension as equation \eqref{equa:indep_2d_BL}. The results are presented in Table \ref{table: 2d_levy_Diag_0110}.
\begin{equation}
    \begin{aligned}
\label{equa:indep_2d_BL}
        d X_{1,t} &= (X_{1,t} - X_{1,t}^3) dt + d B_{t} \\
        d X_{2,t} &= (X_{2,t} - X_{2,t}^3) dt + d L_{t}\,,
    \end{aligned}
\end{equation}
\end{example}

\begin{table}[ht]
\centering 
\scalebox{0.9}{
\begin{tabular}{ccccccc}
\toprule Parameter & $\lambda_{0}^{(1)}$ & $\lambda_{1}^{(1)}$ & $\lambda_{2}^{(1)}$ & $\lambda_{3}^{(1)}$ & $\sigma^{(1)}$ & $\xi^{(1)}$ \\
\hline
True & 0 & 1 & 0 & $-1$ & 1 & 0 \\
\hline
Estimated & $0.0$ & $1.0802537$ & $0.0$ & $-1.0478958$ & $0.9991456$ & $0.0$ \\
\hline
Parameter & $\lambda_{0}^{(2)}$ & $\lambda_{1}^{(2)}$ & $\lambda_{2}^{(2)}$ & $\lambda_{3}^{(2)}$ & $\sigma^{(2)}$ & $\xi^{(2)}$ \\
\hline
True & 0 & 1 & 0 & $-1$ & 0 & 1 \\
\midrule
Estimated & $0.0$ & $1.1367049$ & $0.0$ & $-1.177175$ & $0.0$ &  $1.0660018$ \\
\bottomrule
\end{tabular}
}
\caption{Reveal unknown drift and diffusion terms of equation \eqref{equa:indep_2d_BL} with $10,000$ samples at time snapshots $t = \{0.1, 0.3, 0.5, 0.7, 1.0\}$ in $2d$.} 
\label{table: 2d_levy_Diag_0110}
\end{table}

\subsubsection{Coupled Dimensions}
\label{sec: couple_}

Sombrero potential $\bm{V} (\bm{x}) = -5\|\bm{x}\|^2 +  \|\bm{x}\|^4= -5(x_1^2 + x_2^2) + (x_1^2 + x_2^2)^2$ is a well-known rotation invariant potential that the symmetry breaking is triggered in the quantum mechanics, and its gradient acts as the force of the stochastic system. We assume that the stochastic system follows the SDE with coupled drift terms $\bm{m} (\bm{x}) = - \nabla \bm{V}$,
\begin{equation}
    \begin{aligned}
        \label{equa: dependent_2d}
        d X_{1,t} &= (10X_{1,t}-4X_{1,t}^3-4X_{1,t}X_{2,t}^2) dt + d L_{1t} \\
        d X_{2,t} &= (10X_{2,t}-4X_{2,t}X_{1,t}^2-4X_{2,t}^3) dt + d L_{2t}\,.
    \end{aligned}
\end{equation}

The basis set for mutually dependent drift terms in equation \eqref{equa: dependent_2d} according to our method is $\Lambda= \{1, x_1, x_2, x_1^2, x_1 x_2, x_2^2, x_1^3, x_1 x_2^2, x_1^2 x_2, x_2^3\}$, and for diffusion term is $\{1, 1\}$, thus the objective function for regression is expressed as equation \eqref{equa: dependent2d'}. Based on polynomial basis, our goal is to estimate the unknown parameters $\{\{\lambda_{i}^{(j)}\}_{i=1}^{b}\}_{j=1}^{2}$, and $\xi^{(j)}$, $j =1, 2$. We sample $20,000$ data points at time snapshots $t = 0.0, 0.1, \cdots, 0.9, 1.0$ to implement WCR method. Table \ref{table: 2d_dependent} implies that our method can also efficiently learn coupled drift terms. 
\begin{equation}
    \begin{aligned}
        \label{equa: dependent2d'}
        d X_{j,t} &= \big(\lambda_0^{(j)} + \lambda_1^{(j)} X_{1,t} + \lambda_2^{(j)} X_{2,t} + \lambda_3^{(j)} X_{1,t}^2 + \lambda_4^{(j)} X_{1,t} X_{2,t} + \lambda_5^{(j)} X_{2,t}^2 \\
        &+ \lambda_6^{(j)} X_{1,t}^3 
        + \lambda_7^{(j)} X_{1,t} X_{2,t}^2 + \lambda_8^{(j)} X_{1,t}^2 X_{2,t} + \lambda_9^{(j)} X_{2,t}^3\big) dt + \xi^{(j)} d L_{j, t} \,.
    \end{aligned}
\end{equation}

\begin{table}[htp]
\centering 
\scalebox{0.9}{
\begin{tabular}{cccccccccccc}
\toprule 
Parameter & $\lambda_{0}^{(1)}$ & $\lambda_{1}^{(1)}$ & $\lambda_{2}^{(1)}$ & $\lambda_{3}^{(1)}$ & $\lambda_{4}^{(1)}$ & $\lambda_{5}^{(1)}$ & $\lambda_{6}^{(1)}$ & $\lambda_{7}^{(1)}$ & $\lambda_{8}^{(1)}$ & $\lambda_{9}^{(1)}$ & $\xi^{(1)}$ \\
\hline
True & 0 & 10 & 0 & 0 &0 & 0 & -4 & 0 & -4 &0 &1 \\
\hline
Estimated & $0.0$ & $10.8513$ & $0.0$ & $0.0$ & $0.0$ & $0.0$ & $-4.3692$ & $0.0$ & $-4.3311$ & $0.0$ & $0.992549$ \\
\midrule Parameter & $\lambda_{0}^{(2)}$ & $\lambda_{1}^{(2)}$ & $\lambda_{2}^{(2)}$ & $\lambda_{3}^{(2)}$ & $\lambda_{4}^{(2)}$ & $\lambda_{5}^{(2)}$ & $\lambda_{6}^{(2)}$ & $\lambda_{7}^{(2)}$ & $\lambda_{8}^{(2)}$ & $\lambda_{9}^{(2)}$ & $\xi^{(2)}$ \\
\hline
True & 0 & 0 & 10 & 0 & 0 & 0 & 0 & -4 & 0 & -4 & 1 \\
\hline
Estimated & $0.0$ & $0.0$ & $11.2223$ & $0.0$ & $0.0$ & $0.0$ & $0.0$ & $-4.4875$ & $0.0$ & $-4.5451$ & $0.9211781$ \\
\bottomrule
\end{tabular}
}
\caption{Estimation of the coupled system with L\'{e}vy noise with $20,000$ samples in $2d$ case. The ground truth of SDE is \eqref{equa: dependent_2d}; The Gaussian number $M =100$.} 
\label{table: 2d_dependent}
\end{table}

\begin{figure}[ht]
\centering    
\includegraphics[width=4.8cm,height=3.6cm]{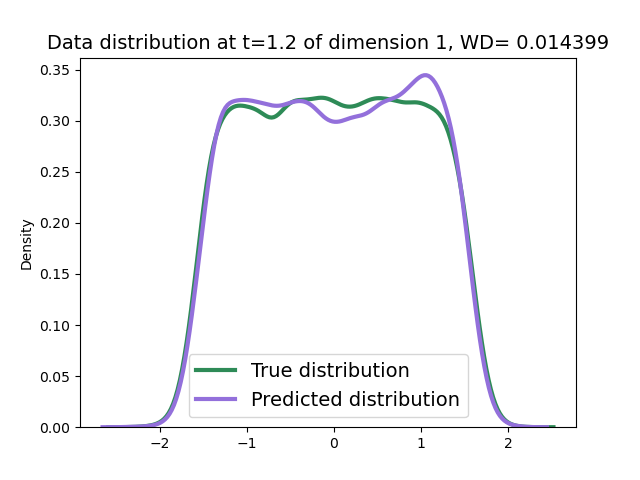}
\includegraphics[width=4.8cm,height=3.6cm]{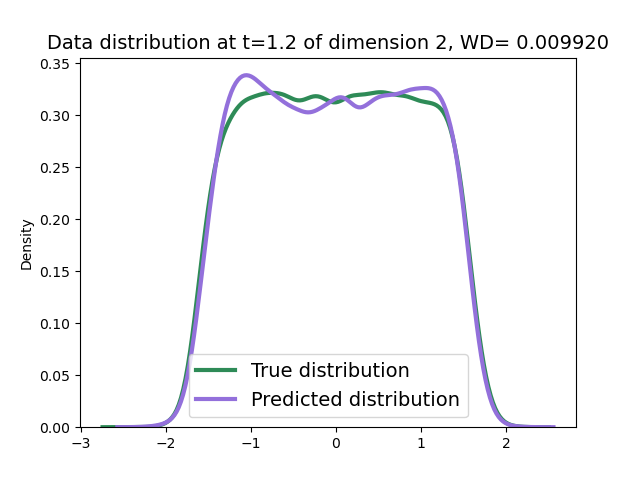}
\includegraphics[width=4.8cm,height=3.6cm]{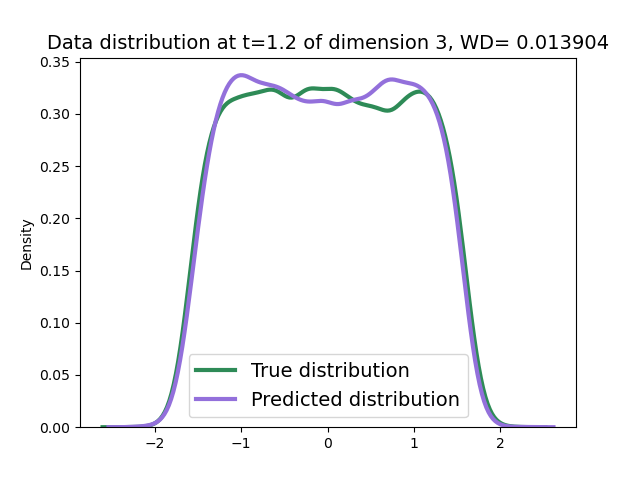}

\caption{Comparison of the true data distribution and the predicted distribution of the estimated SDE at the end moment $t=1.2$ for $3d$ coupled system. The "WD" in figures means "Wasserstein Distance" of the predicted one and the ground truth.}
\label{fig: dist_couple3d}
\end{figure}

Similarly, we consider the potential $\bm{V}$ in three dimension. We choose basis set as $\Lambda = \{ 1, x_1, x_2, x_3, x_1^2, x_1 x_2, x_1 x_3, x_2^2, x_2 x_3, x_3^2, x_1^3, \cdots, x_3^3 \}$, which contains $20$ elements, and Figure \ref{fig: dist_couple3d} gives the comparison at $t=1.2$ of true data distribution and the predicted distribution of estimated SDE. The predicted distribution is the same as Section \ref{sec: necessity}. The "WD" means "Wasserstein-1 Distance" between two sets of discrete data points, calculated by python package. The distribution plots are generated by kernel density estimation in Python, and the MRE for total $63$ parameters is $0.13865$.

\subsection{High Dimensional Problems}
\label{sec:5d}
We also conduct experiments in five dimensions, while previous studies typically only perform two-dimensional experiments with L\'{e}vy noise. Utilizing the confluent hypergeometric function \cite{burkardt2021unified}, the integral term can be solved easily in WCR method. Consider the diagonal matrix $\bm{\xi} (\bm{X_t}) = diag\{ \xi_1, \xi_2, \cdots, \xi_d
 \}$ and each $\xi_i$ are constants. We separately test the independent case and coupled system. In both scenarios, we used a sample size of $30,000$, with training samples being uniformly observed at $t=0, 0.1, \cdots, 1.0$.
\begin{equation}
    \begin{aligned}
    \label{equa:indep_dim5}
    d \bm{X_t} =
    d \left(
        X_{1t},\cdots, X_{dt}
        \right)^T= \bm{m}( \bm{X_t}) dt + d \left(
        B_{1t},\cdots, B_{dt}
        \right)^T + \text{diag}\{\xi_1,\cdots, \xi_d\}
        d\left(
        L_{1t},\cdots, L_{dt}
        \right)^T
    \end{aligned}
\end{equation}

\paragraph{Case One: $\bm{m}(\bm{X}_t) = (X_{1,t} - X_{1,t}^3, \cdots, X_{5,t} - X_{5,t}^3)^T$, $\bm{\sigma} = \bm{\xi} = diag \{1, \cdots, 1\}_{5\times 5}$.} 
Our approach does not need a sample size increasing with the geometric growth of dimensions, even in high-dimensional independent scenarios, unlike numerical methods. Table \ref{table: 5d_levy} demonstrates positive results.

\begin{table}[htp]
\centering 
\scalebox{0.9}{
\begin{tabular}{ccccccc}
\toprule 
True Parameter & 0 & 1 & 0 & $-1$ & 1 & 1 \\
\hline
Parameter & $\lambda_{0}^{(1)}$ & $\lambda_{1}^{(1)}$ & $\lambda_{2}^{(1)}$ & $\lambda_{3}^{(1)}$ & $\sigma^{(1)}$ & $\xi^{(1)}$ \\
\hline
Estimated & $-0.12000132$ & $0.93020815$ & $-0.10496778$ & $-1.1212871$ & $1.0527982$ & $1.0197511$ \\
\hline
Parameter & $\lambda_{0}^{(2)}$ & $\lambda_{1}^{(2)}$ & $\lambda_{2}^{(2)}$ & $\lambda_{3}^{(2)}$ & $\sigma^{(2)}$ & $\xi^{(2)}$ \\
\hline
Estimated & $-0.07088655$ & $1.0873672$ & $-0.06265919$ & $-1.1261857$ & $1.1448588$ & $0.9033003$ \\
\hline
Parameter & $\lambda_{0}^{(3)}$ & $\lambda_{1}^{(3)}$ & $\lambda_{2}^{(3)}$ & $\lambda_{3}^{(3)}$ & $\sigma^{(3)}$ & $\xi^{(3)}$\\
\hline
Estimated & $-0.05240659$ & $0.999507$ & $0.06522602$ & $-0.9960823$ & $1.1572613$ & $0.90616125$ \\
\midrule Parameter & $\lambda_{0}^{(4)}$ & $\lambda_{1}^{(4)}$ & $\lambda_{2}^{(4)}$ & $\lambda_{3}^{(4)}$ & $\sigma^{(4)}$ & $\xi^{(4)}$\\
\hline
Estimated & $-0.09751202$ & $1.041346$ & $-0.10804728$ & $-1.1485909$ & $1.1243781$ & $0.94195074$ \\
\hline
Parameter & $\lambda_{0}^{(5)}$ & $\lambda_{1}^{(5)}$ & $\lambda_{2}^{(5)}$ & $\lambda_{3}^{(5)}$ & $\sigma^{(5)}$ & $\xi^{(5)}$\\
\hline
Estimated & $-0.04412302$ & $1.1565177$ & $-0.08664627$ & $-1.1861836$ & $1.1126705$ & $0.9244608$ \\
\bottomrule
\end{tabular}
}
\caption{Reveal the unknown parameters in $5d$ independent problems. The Gaussian number is set as $M =600$. The training time on $CPU$ device is around $700 s$.} 
\label{table: 5d_levy}
\end{table}



\begin{figure}[h]
\centering    
\includegraphics[width=4.85cm,height=3.6cm]{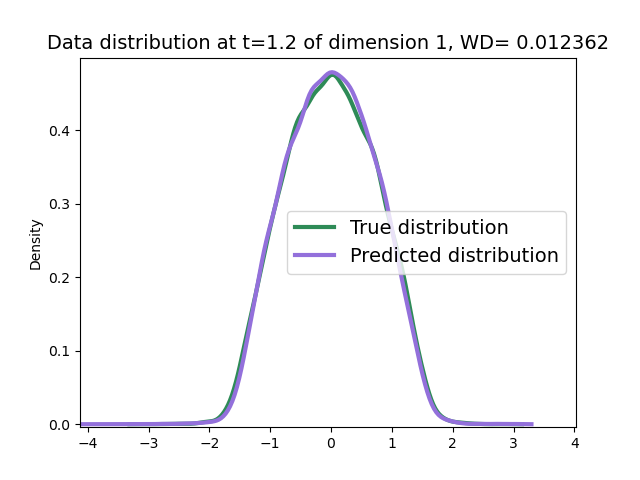}
\includegraphics[width=4.85cm,height=3.6cm]{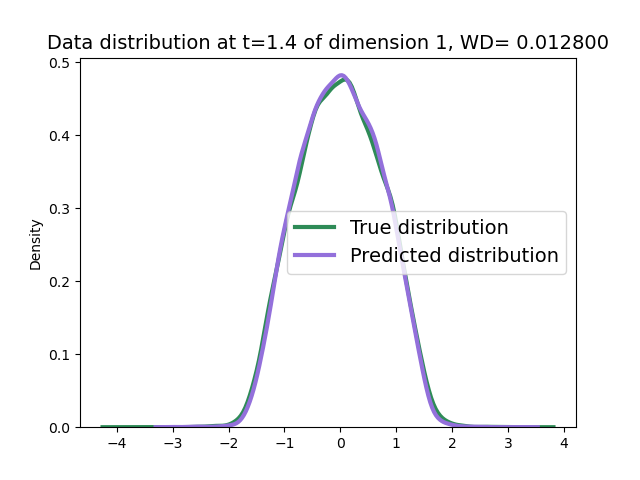}
\includegraphics[width=4.85cm,height=3.6cm]{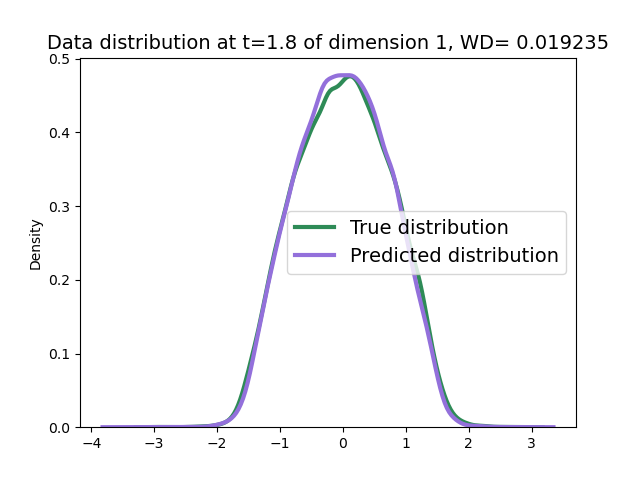}
\caption{Comparison of the true data distribution and the predicted distribution of the estimated SDE at moments $t=1.2, 1.4, 1.8$ for $5d$ coupled system. The "WD" in sub-titles means "Wasserstein Distance" of the predicted one and the ground truth.}
\label{fig: dist_couple_5d}
\end{figure}

\paragraph{Case Two: $\bm{m}(\bm{X}_t) = -\nabla \bm{V}(\bm{X}_t)$, where $\bm{V}(\bm{X}_t) = -5 \|\bm{X}_t \|^2 + \|\bm{X}_t \|^4$, $\bm{\sigma} = \bm{0}$, $\bm{\xi} = diag\{1, \cdots,1\}_{5\times 5}$.} In constructing the basis for the five-dimensional coupled system, we simplify the basis as $\{1, x_1, \cdots, x_5, x_1^3, x_1 x_2^2, \cdots, x_5^3 \}$, which contains $31$ terms, as the drift term only consists of first-order and third-order polynomials. Note that with the increase in dimension, the number of polynomial basis terms escalates significantly, leading to an exponential rise in the size of the sample space. In our experiment, we selected $4,000$ Gaussian sampling points. Since we consider systems with rotational invariance, we plot the results in one dimension as an example. The unknown parameters in the SDE are estimated by WCR given snapshots at $t = \{0, 0.1, \cdots, 1.0\}$, and then we use the plots of the distributions of the real data for $t=1.2, 1.4, 1.8$ as a measure of the predicted outcome of the estimated SDE. Results in Figure \ref{fig: dist_couple_5d} suggest that our approach achieves good estimation for prediction. We can also see that in the presence of increases in dimension, WCR still yields accurate estimates even when the sample size $N$ does not increase significantly, remaining comparable to the sample size in $2d$ and $3d$.

\subsection{Non-polynomial Drift Terms}
\label{sec: nonpoly}

The SDE with a periodic drift in a velocity field has some applications \cite{bhat2022drift}. To expand the basis dictionary, we introduce sinx and cosx into the basis, and consider the SDE $dX_t = (-X_t - \sin X_t) dt + dB_t + dL_t$. The basis is selected as $\{1, x, x^2, x^3, \sin x, \cos x, \sin^3 x, \cos^3 x\}$. Parameters are estimated from the SDE $dX_t = \lambda_{0} + \lambda_{1} X_t + \lambda_{2} X_t^2 + \lambda_{3} X_t^3 + \lambda_{4} \sin x + \lambda_{5} \cos x + \lambda_{6} \sin^3 x + \lambda_{7} \cos^3 x + \sigma dB_t + \xi dL_t$. The results of the improved dictionaries are shown in Table \ref{table: sinDrift}, highlighting the widespread application properties of the WCR.

\begin{table}[htpb]
\centering
\scalebox{0.9}{
\begin{tabular}{ccccccccccc}
\toprule Parameter & $\lambda_{0}$ & $\lambda_{1}$ & $\lambda_{2}$ & $\lambda_{3}$ & $\lambda_{4}$ & $\lambda_{5}$ & $\lambda_{6}$ & $\lambda_{7}$ & $\sigma$ & $\xi$ \\
\hline
True & 0 & $-1$ & 0 & 0 & $-1$ & 0 & 0 & 0 & 1 & 1 \\
\hline
Estimated & $0.0$ & $-0.961751$ & $0.0$ & $0.0$ & $-1.025329$ & $0.0$ & $0.0$ & $0.0$ & $0.9720182$ & $1.0066$ \\
\bottomrule
\end{tabular}
}
\caption{Estimated results for non-polynomial drift term $m(x) = -x- \sin x$, the diffusion term $\sigma$, and the noisy intensity $\xi$. The sample size is $20,000$.} 
\label{table: sinDrift}
\end{table}

\begin{figure}[htbp]
\label{fig:non-poly}
\centering
\includegraphics[width=6cm,height=4.8cm]{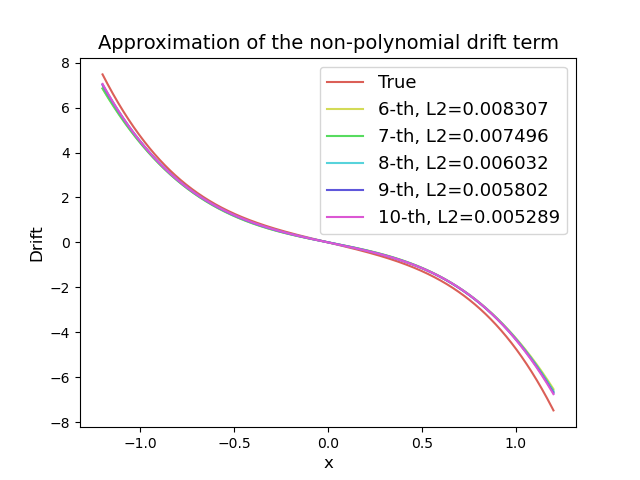}
\caption{The estimated drift terms of SDE with non-polynomial drift term $m(x) = - 4x^3-2x e^{-x^2}$. 
}
\end{figure}

In cases where the specific form of the drift term is unknown, the WCR approach can still be utilized to infer hidden laws. By selecting polynomial basis functions, we can approximate non-polynomial drift terms. For instance, the SDE $dX_t = (-4X_t^3 - 2X_t e^{-X_t^2}) dt + dL_t$ is considered, and the WCR is employed to obtain an approximating polynomial drift term to replace the original drift. A data sample of $20,000$ points for each snapshot is taken at time points $t=0.1, \cdots, 1.0$. The estimated results in Figure \ref{fig:non-poly} illustrate that the approximation errors consistently decrease from the sixth to the tenth order polynomials. To quantify the approximation error, the relative $2$-norm error, denoted as $L_2$ error, is calculated by $L_2 = \frac{\|\hat{m}(\bm{x}) - m(\bm{x}) \|_2}{\|m(\bm{x}) \|_2}$, where $\hat{m}$ represents the estimated polynomial drift term and $\bm{x}$ denotes random sample tensors within the bounded area $[-1.2, 1.2]$. The $L_2$ error of approximation in $[-1.2, 1.2]$ are shown in legends. These two examples collectively demonstrate the applicability of the WCR method.

\subsection{Unconstant Diffusion Terms}
\label{sec:func_diff}

Above experiments consider the constants $\sigma$ and $\xi$. Black-Scholes (BS) model is one of the most classic models for valuation of options and hedging \cite{black1973pricing}, and it assumes that the underlying stock price follows the Geometric Brownian motion (GBM) $d X_t = m (X_t) dt + \sigma (X_t) d B_t $, in which the diffusion term $\sigma(X_t)$ is a linear function of $X_t$. We then illustrate that the WCR can also reveal stochastic dynamics with noise intensity dependent on $X_t$. In this section, $10,000$ data points are sampled at non equidistant moments $t = 0.1, 0.3, 0.4, 0.7, 1.0$.

\begin{example}[Geometric Brownian motion] Consider
\begin{equation}
    \begin{aligned}
    \label{eq: BS_gaussian}
        d X_t = m_1(X_t) dt + \sigma \cdot X_t d B_t \,,
    \end{aligned}
\end{equation}
where $m_1 (X_t), \sigma \in \mathbb{R}$, $B_t$ is a standard Brownian motion. The Fokker-Planck equation for equation \eqref{eq: BS_gaussian} is $\frac{\partial}{\partial t} p(x, t) 
    = -\frac{\partial}{\partial x} \left(m_1(x) p(x, t)\right) + \frac{\partial^2}{\partial x^2} \left(\frac{1}{2} \sigma^2 x^2 p(x, t) \right)\,.$
\end{example}
Choose basis $\Lambda^m = \{ 1, x, x^2, x^3\}, \Lambda^\sigma = \{ 1, x, x^2\}$ for $m_1$ and $G = \sigma^2 x^2$ respectively, and we expand $m$ and $G$ as a linear combination of basis, 
$$m_1 = \lambda_0^m + \lambda_1^m x +\lambda_2^m x^2 + \lambda_3^m x^3\,, \quad G = \lambda_0^\sigma + \lambda_1^\sigma x + \lambda_2^\sigma x^2 \,,$$
we then estimate $\{\lambda_i^m \}, \{ \lambda_j^\sigma\}, (i = 0, 1, 2, 3; j=0, 1, 2)$. Table \ref{table: BS_GBM} shows the applicability of WCR algorithm to the Black-Scholes model, which provides ideas for sparse regression methods to quickly solve some financial practical problems.

\begin{table}[htp]
\centering 
\scalebox{0.9}{
\begin{tabular}{cccccccc}
\toprule Parameter & $\lambda_{0}^m$ & $\lambda_{1}^m$ & $\lambda_{2}^m$ & $\lambda_{3}^m$ & $\lambda_0^\sigma$ & $\lambda_1^\sigma$ & $\lambda_2^{\sigma}$ \\
\hline
True parameters & 0 & 1 & 0 & $-1$ & 0 & 0 & 1\\
\hline
Estimated parameters & $0.0$ & $0.9674$ & $0.0$ & $-0.9115$ & $0.0$ & $0.0$ & $0.9158$ \\
\bottomrule
\end{tabular}
}
\caption{Estimation of SDE with linear functional diffusion term. The ground truth is $d X_t = (X_t -X_t^3) dt + X_t \cdot d B_t$.} 
\label{table: BS_GBM}
\end{table}

\begin{example}
GBM, however, does not describe well extreme events, such as dramatic price drops, then many generalizations of BS models that incorporate $\alpha$-stable L\'{e}vy noise have been introduced \cite{aguilar2017non}. 
We then test our method on equation \eqref{eq: BS_levy}.
\begin{equation}
    \begin{aligned}
    \label{eq: BS_levy}
        d X_t = m_2 (X_t) dt + \xi \cdot X_t d L_t \,,
    \end{aligned}
\end{equation}
where $m_2 (X_t), \xi \in \mathbb{R}$, $L_t$ is an $\alpha$-stable symmetric L\'{e}vy process. The Fokker-Planck equation for equation \eqref{eq: BS_levy} is 
$\frac{\partial}{\partial t} p(x, t) 
    = -\nabla\cdot \left( m_2(x) p(x, t)\right) - |\xi|^{\alpha} (-\Delta)^{\alpha/2} ( |x|^{\alpha} p(x, t))\,.$
\end{example}
We choose $\Lambda^m = \{ 1, x, x^2, x^3\}, \Lambda^\xi = \{ 1, |x|^\alpha \}$, and expand $m(x)$ and $|\xi \cdot x|^{\alpha}$ as a linear combination of basis, 
    \begin{align*}
        m(x) = \lambda_0^m + \lambda_1^m x +\lambda_2^m x^2 + \lambda_3^m x^3\,, \quad |\xi \cdot x|^{\alpha} = \lambda_0^\xi + \lambda_1^\xi |x|^{\alpha} \,,
    \end{align*}
then estimate $\{\lambda_i^m \}, \{ \lambda_j^\xi\}, (i = 0, 1, 2, 3; j=0, 1)$. The Gaussian number is set as $M =70$, more than GBM, since L\'{e}vy noise leads to a heavy-tailed distribution. Specific results are shown in Table \ref{table: BS_GBM_levy}. Here we have only done the non-constant noise intensity problem in one-dimensional case, since the FP equation has the term $|\bm{\xi}(\bm{X}_t)|^\alpha$, difficult to disentangle as a polynomial basis in high dimension.

\begin{table}[htp]
\centering 
\scalebox{0.9}{
\begin{tabular}{ccccccc}
\toprule Parameter & $\lambda_{0}^m$ & $\lambda_{1}^m$ & $\lambda_{2}^m$ & $\lambda_{3}^m$ & $\lambda_0^\xi$ & $\lambda_1^\xi$ \\
\hline
True parameters & 0 & 1 & 0 & $-1$ & 0 & 1 \\
\hline
Estimated parameters & $0.0$ & $0.9080506$ & $0.0$ & $-1.039701$ & $0.0$ & $0.9089$ \\
\bottomrule
\end{tabular}
}
\caption{Reveal the unknown terms for SDE with functional noise intensity $L_t$. The ground truth is $d X_t = (X_t -X_t^3) dt + X_t d L_t$.} 
\label{table: BS_GBM_levy}
\end{table}




\subsection{Robustness}
\label{sec:robust}
We test the robustness of our approach by conducting experiments with noisy data. Two types of random noise, addictive and multiplicative noises, with scale $\delta$, are added respectively. 
\begin{equation}
\begin{aligned}
\hat{\bm{x}}_{t_\ell}^k = \bm{x}_{t_\ell}^k + \delta\mathcal{U}_{t_\ell}^k\,, \quad \text{and} \,\,
\hat{\bm{x}}_{t_\ell}^k = \bm{x}_{t_\ell}^k (1+ \delta\mathcal{U}_{t_\ell}^k) \,,
\end{aligned}
\end{equation}
where $\mathcal{U}_{t_\ell}^k$ is an uniform random variable in $[-1, 1]$, $\ell =\{1, \cdots,L\}$, and $k = 1, \cdots, N$. 
We sample $N=10,000$ data points in each snapshot with time point set $\{0, 0.1, \cdots, 1.0 \}$. The MRE of unkbown terms, including the drift, diffusion and noise intensity terms, are showed in Table \ref{table: noisy}. Results imply that our approach has robustness. 

\begin{table}[htp]
    \centering
    \scalebox{0.9}{
    \begin{tabular}{cccccc}
         \toprule
         Additive Noise & 
         $\delta=0\%$ & $\delta = 5\%$ & $\delta=10\%$ & $\delta=20\%$\\
         \midrule
         MRE in drift& $0.0319$
         &  $0.0117$ & $0.0447$ & $0.1041$  \\
         MRE in diffusion & $0.0004$ & $0.0110$& $0.0261$ & $0.0380$   \\
         \midrule
         Multiplicative Noise &
         $\delta=0\%$ & $\delta = 5\%$ &
         $\delta=10\%$ & $\delta=20\%$ \\
         \midrule
         MRE in drift & $0.0319$ & $0.0514$ & $0.1056$ & $0.2957$  \\
         MRE in diffusion & $0.0004$ & $0.0345$ & $0.0178$ & $0.0608$   \\
         \bottomrule
    \end{tabular}
    }
    \caption{The estimated results of one-dimensional problem with noisy data. Set $\delta=5,10,20\%$. The ground truth is $d X_t= (X_t-X_t^3) dt + dB_t +dL_t$.}
    \label{table: noisy}
\end{table}


\section{Conclusions}
\label{sec:conclu}
This study aims to develop a method for learning the stochastic dynamics that involve both Gaussian and L\'{e}vy noises, only from the aggregate data. The inclusion of L\'{e}vy noise in the SDE leads to the presence of a computationally demanding integral term in the Fokker-Planck equation. Moreover, the simultaneous learning of these two noise types poses challenges, particularly the impediment to differentiation arising from their interaction. Notably, our approach effectively solves the aforementioned challenges, as evidenced by the experimental results. Furthermore, compared to prior works \cite{chen2021solving,gao2016fokker,li2022extracting} that conduct at most two-dimensional experiments with L\'{e}vy noise, our method exhibits effectiveness in solving a five-dimensional coupled system.

Our experiments achieve less than $1.0 s$ in one-dimensional problems and around $10 s$ in two dimension on $CPU$ device. As the dimension increases, the required Gaussian kernel functions increase and the dimensionality of the coefficient matrix for the regression rises, leading to more computational burdens. Notably, WCR does not experience a geometric increase in sampling points as the problem dimension increases. The rapid recognition capability of our approach can be applied in the pre-training of other models.

\subsection{Discussions}

\paragraph{Source of error.}
Since our approach includes Monte Carlo and temporal differencing, increasing the number of samples and snapshots naturally improves the accuracy. In addition, the effect of Gaussian sampling points cannot be ignored. The Gaussian number $M$ affects the results since the regression coefficient matrix $\bm{A}$ is in $\mathbb{R}^{LM\times b}$. Large $M$ results in a more sparse coefficient matrix and an unstable solution. Conversely, if $M$ is assigned a smaller value, Gaussian samples inadequately cover the data points, leading to information gaps and sub-optimal estimation.  Specific error analysis of collocation points could be similar to \cite{lu2024weak,messenger2022learning}. Some other factors such as learning two noises at the same time will be more difficult than learning only one noise due to the two noises interacting with each other, and restricting the L\'{e}vy variables to a bounded area causes missing information.

\paragraph{Sampling.} 
The error in high dimension is larger than low dimensions. This is primarily due to the exponential increase in the number of basis $b$ for higher-dimensional cases, subsequently leading to rapid dimensional expansion of the coefficient matrix $\mathbb{R}^{LM\times b}$ and an increase in matrix sparsity. Furthermore, since we sample Gaussian functions based on the data distribution in the full space, higher dimensions necessitate a greater samplings. If an adaptive sampling approach is explored, we do not need to sample so many Gaussians in full space instead of less samplings at each moment. It speeds up the algorithm and possibly reduces the error. 
On the other hand, the sampled kernels remain fixed or static, whereas the data points continue to evolve over time. This dynamic nature of data points suggests that a dynamic sampling scheme could provide monitoring of all evolving points, thereby enhancing error reduction. As a future consideration, an adaptive sampling approach could be explored to reduce computational complexity and improve accuracy in high-dimensional systems \cite{zhao2022adaptive}.


\section*{Acknowledgements}
This work was supported by the National Key R\&D Program of China (Grant No. 2021YFA0719200) . The authors would like to thank the helpful discussions from Yan Jiang.

\bibliographystyle{plain}
\bibliography{ref}

\appendix
\section{Experimental Details}
\subsection{Data Generating}
\label{appen: data}

\paragraph{Generate stable L\'{e}vy variables.}
We use Euler scheme to generate L\'{e}vy random variables \cite{yang2020generative}.
If $\bm{X}_t \in \mathbb{R}^d$ satisfies SDE (equation \eqref{equa:sde_mf_bmlevy}) with both L\'{e}vy noise and Gaussian noise, the Euler scheme provides that 
\begin{equation}
    \begin{aligned}
        \label{equa: euler_levy}
        & \bm{X}_0 = \mathcal{N}(0, z) \\
        &\bm{X}_{(\ell+1)\Delta t} = \bm{X}_{\ell\Delta t} + \bm{m} \Delta t + \bm{\sigma} \sqrt{\Delta t} \eta_\ell + \bm{\xi} (\Delta t)^{1/\alpha} \zeta_{\alpha, \ell}\,,
    \end{aligned}
\end{equation}
with initial $\bm{X}_0$, where $\eta_\ell$ is randomly sampled from a standard Gaussian distribution, and $\zeta_{\alpha, i}$ are i.i.d. $\alpha$-stable random variables ($\alpha \in (0, 1) \cup (1, 2]$). The method for generating $\alpha$-stable L\'{e}vy random variable $X \sim S_{\alpha} (1, 0, 0)$ is as follows.

(i) Generate a random variable $V$ satisfying uniform distribution $U(-\frac{\pi}{2}, \frac{\pi}{2})$;
    
(ii) Generate another random variable $S$ satisfying uniform distribution $U(0, 1)$, and let $W = -\ln S$ be an exponential random variable independent of $V$ with mean $1$;
    
(iii) Let $X = \frac{\sin(\alpha V)}{(\cos (V))^{1/\alpha}} \times (\frac{\cos(V - \alpha V)}{W})^{(1-\alpha)/\alpha}$\,.

\paragraph{Limit L\'{e}vy variables in a bounded area.} The unbounded computational domain, resulting from the heavy-tailed property and jumps of the L\'{e}vy process, poses significant challenges to algorithm operation. L\'{e}vy noise also has a notable impact on differentiating Gaussian noise in locations that are relatively distant \cite{mantegna1996turbulence}. In order to address these issues, this study introduces small perturbations to the process of generating L\'{e}vy variables. It is worth noting that stable variables are generated by the expression $X = \frac{\sin(\alpha V)}{(\cos (V))^{1/\alpha}} \times (\frac{\cos(V - \alpha V)}{W})^{(1-\alpha)/\alpha}$, as mentioned earlier. To maintain the boundedness of $X$, we introduce noise perturbations to the variable $V$, ensuring that it does not approach $-\pi/2$ and $\pi/2$, thereby preventing the resulting $X$ from reaching infinity.

Some works also propose methods to generate bounded L\'{e}vy data. For instance, some focus on cutting off L\'{e}vy data, such as using the truncated L\'{e}vy distribution (TLD) \cite{matacz2000financial} and directly truncating the data, which involves ignoring data outside the bounded region \cite{yang2020generative,gao2016dynamical}. Another technique is performing function transformations on the data, such as using $\tanh(\cdot)$, to limit it to a bounded range \cite{yang2020generative}. Additionally, one can add noise to the generated data \cite{yang2020generative}.

\subsection{Gaussian Kernels with Different Variances}
\label{appen: improve}

In our study, we begin by extracting $M$ kernels from a $d$-dimensional space using the Latin Hypercube Sampling (LHS) method. This involved the following steps: first, dividing each dimension into non-overlapping intervals, resulting in a total of $m$ intervals, each with an equal probability and length, typically achieved through the use of a uniform distribution. Next, a point is randomly selected from each interval in every dimension, with the selected points combined to form a vector. We also restrict the range of LHS to $l_b + lhs.ratio \times (u_b - l_b)$, where $lhs.ratio$ denotes the proportion of the sampling within the entire range of data points, and $u_b$ and $l_b$ represent the maximum and minimum values of the data points, respectively. Throughout our experiments, we observe that different Gaussian variances have impacts on the estimation. It became evident that different variances should be chosen for different settings.

\begin{table}[htp]
\centering 
\scalebox{0.9}{
\begin{tabular}{ccccccc}
\toprule Parameter & $\lambda_{0}$ & $\lambda_{1}$ & $\lambda_{2}$ & $\lambda_{3}$ & $\sigma$ & $\xi$ \\
\hline
True parameter & 0 & 1 & 0 & $-1$ & 1 & 1 \\
\hline
(a) & $-0.03502786$ & $0.9731436$ & $0.02529686$ & $-1.0175471$ & $0.9148$ & $1.0156$ \\
\hline
(b) & $0.0$ & $0.96525615$ & $0.0$ & $-1.0092138$ & $0.9794$ & $1.0372$ \\
\bottomrule
\end{tabular}
}
\caption{Estimated results where Gaussian kernels have different variance in one dimension. (a) $std_1 = 1, std_2 = 0.4$,
$lhs.rat_1 = 1, lhs.rat_2 = 1.25 $, ratio of Gaussian sample number is $0.7:0.3$; (b) $std_1 = 1, std_2 = 0.5$, $lhs.rat_1 = 1.5, lhs.rat_2 = 0.9$, ratio of Gaussian sample number is $ 0.7:0.3$; We set $20,000$ samples, observations are taken at $t = 0.1, 0.2, \cdots, 0.9, 1.0$. The samples adhere the SDE $d X_t= (X_t-X_t^3) dt + dB_t + dL_t $.} 
\label{table: both_1d}
\end{table}

In testing the WCR, we choose Gaussian functions with two different variances. The Gaussian mean values are selected based on the data distribution, and the variances are set as $1$ and $0.4$, with corresponding $lhs.ratio$ values of $1$ and $1.25$, respectively. 
The result below leads us to consider the potential of adaptively or dynamically sampling Gaussian kernels, further underscoring the significance of the experimentation process.


\section{Weak Form of the Fokker-Planck Equation}
\label{appen: weak}
We consider a more general case that the noise intensity for rotational invariant L\'{e}vy noise $\bm{\xi}(\bm{x})$ ia a diagonal matrix with scalar function diagonal elements rather than constants. Recall that the FP equation derived from a SDE (Remark \ref{remark:RI}) is
\begin{equation}
    \begin{aligned}
    \label{equa:fp}
    \frac{\partial}{\partial t} p(\bm{x}, t) = -div \left(\bm{m}(\bm{x}) p(\bm{x}, t)\right) + \sum_{i,j}^{d} \partial_{ij} \left( G_{i,j}(\bm{x}) p(\bm{x}, t) \right) + \theta_{\alpha, d} (-\Delta)^{\alpha/2} \left( |\xi_0(\bm{x})|^{\alpha/2} p(\bm{x}, t)\right)\,.
\end{aligned}
\end{equation}

If the fractional derivative term is not included, we can establish the weak form through integration by parts.
For SDE driven by L\'{e}vy process, $d \bm{X}_t = \bm{m}(\bm{X}_t) dt + \bm{\xi}(\bm{X}_t) d \bm{L}_t$ and its FP equation 
$\partial_t p(\bm{x},t) + div (\bm{m}(\bm{x}) p(\bm{x},t)) + \theta(-\Delta)^{\alpha/2}(|\bm{\xi}(\bm{x})|^{\alpha} p(\bm{x},t)) = 0$, Wei et.al. \cite{wei2015well} prove the existence of the weak solution of this PIDE, and we give a simple proof from Fourier transform. Since $\psi(\bm{x})$ is a Gaussian function, it satisfies that $\psi(\bm{x}) = \psi(2\rho -\bm{x})$, where $\rho$ is the mean value of this Gaussian function. Then

\begin{equation}
\begin{aligned}
& \int_{\mathbb{R}^d} \psi(\bm{x}) \cdot (-\Delta)^{\alpha/2} \left( |\xi_0(\bm{x})|^\alpha p(\bm{x})\right) d \bm{x}= \int_{\mathbb{R}^d} \psi(2\rho -\bm{x}) \cdot (-\Delta)^{\alpha/2} \left( |\xi_0(\bm{x})|^\alpha p(\bm{x})\right) d \bm{x} \\
&= (-\Delta)^{\alpha/2} \left( |\xi_0(2\rho)|^\alpha p(2\rho)\right) \ast \psi(2\rho) =  \cF^{-1} \left[ \widehat{\bigg((-\Delta)^{\alpha/2} (|\xi_0|^\alpha p)\bigg)} \cdot \hat{\psi} \right] (2\rho)\\
&= \cF^{-1} \left[ \hat{\psi} \cdot \left(-\|\bm{k}\|^\alpha \widehat{ \bigg(|\xi_0(\bm{k})|^\alpha p(\bm{k}) \bigg)}\right)\right](2\rho) =
\cF^{-1} \left[ \widehat{ \bigg(|\xi_0(\bm{k})|^\alpha p(\bm{k}) \bigg)}\cdot \left(-\|\bm{k}\|^\alpha \hat{\psi}\right)\right](2\rho)\\
&= \cF^{-1} \left[ \widehat{ \bigg(|\xi_0(\bm{k})|^\alpha p(\bm{k}) \bigg)}\cdot \widehat{(-\Delta)^{\alpha/2} \psi(2\rho)} \right] = \bigg(|\xi_0(2\rho)|^\alpha p(2\rho) \bigg) \ast (-\Delta)^{\alpha/2} \psi(2\rho) \\
&= \int_{\mathbb{R}^d} \left( |\xi_0(\bm{x})|^\alpha p(\bm{x})\right) \cdot (-\Delta)^{\alpha/2} \psi(\bm{x}) d \bm{x}\,.
    \end{aligned}
\end{equation}

Then integration by parts also holds for fractional Laplace operators. For independent L\'{e}vy noise and constant diagonal elements, the weak formula is naturally satisfied.
Therefore, by the independence of $\bm{L}_t$ and $\bm{B}_t$, the weak form of Fokker-Planck equation derived by SDE with both Gaussian noise and L\'{e}vy noise is obtained. 
\begin{equation}
\begin{aligned}
\frac{d}{d t} \int_{\mathbb{R}^d} p(\bm{x}, t) \psi(\bm{x}) d\bm{x}
    &= \int_{\mathbb{R}^d} \bm{m}(\bm{x}) p(\bm{x}, t) \nabla \psi(\bm{x}) d\bm{x} + \sum_{i,j}^{d} \int_{\mathbb{R}^d} p(\bm{x}, t) G_{i,j}(\bm{x}) \partial_{ij}\psi(\bm{x}) d\bm{x}\\
    & - \theta_{\alpha, d} \int_{\mathbb{R}^d} (-\Delta)^{\alpha/2} \psi(\bm{x}) \cdot |\xi_0(\bm{x})|^{\alpha}  p(\bm{x}, t) d\bm{x}\,.
\end{aligned}
\end{equation}

\end{document}